\documentclass{elsart}

\usepackage{graphics}
\usepackage{graphicx}
\usepackage{amssymb}
\usepackage{amsmath}
\usepackage{xcolor}
\usepackage[all]{xy}
\begin{document}

\begin{frontmatter}

% Title, authors and addresses

% use the thanksref command within \title, \author or \address for footnotes;
% use the corauthref command within \author for corresponding author footnotes;
% use the ead command for the email address,
% and the form \ead[url] for the home page:
% \title{Title\thanksref{label1}}
% \thanks[label1]{}
% \author{Name\corauthref{cor1}\thanksref{label2}}
% \ead{email address}
% \ead[url]{home page}
% \thanks[label2]{}
% \corauth[cor1]{}
% \address{Address\thanksref{label3}}
% \thanks[label3]{}

\title{Mathematical analysis and numerical resolution of a 
heat transfer problem arising in water recirculation%\thanksref{label1}
}

%\thanks[label1]{Supported by project MTM2015-65570-P of MINECO/FEDER (Spain).}

% use optional labels to link authors explicitly to addresses:
% \author[label1,label2]{}
% \address[label1]{}
% \address[label2]{}

\author[Santiago]{Francisco J. Fern\'andez},
\ead{fjavier.fernandez@usc.es}
\author[Vigo]{Lino J. Alvarez-V\'azquez\corauthref{cor}},
\corauth[cor]{Corresponding author. Tel: +34 986 812166. Fax: +34 986 812116.}
\ead{lino@dma.uvigo.es}
\author[Vigo]{Aurea Mart\'{\i}nez}
\ead{aurea@dma.uvigo.es}

\address[Santiago]{Universidade de Santiago de Compostela, Instituto de 
Matem\'aticas, 15782 Santiago, Spain.}
\address[Vigo]{Universidade de Vigo, E.I. Telecomunicaci\'on, 36310 Vigo, Spain.}

\begin{abstract}
This work is devoted to the analysis and resolution of a well-posed mathematical model for several processes 
involved in the artificial circulation of water in a large waterbody. 
This novel formulation couples the convective heat transfer equation with the modified Navier-Stokes system following a Smagorinsky turbulence model, 
completed with a suitable set of mixed, nonhomogeneous boundary conditions of diffusive, convective and radiative type.
We prove several theoretical results related to existence of solution, and propose a full algorithm for its computation,
illustrated with some realistic numerical examples.
\end{abstract}

\end{frontmatter}

%\newpage

\section{Introduction}

Artificial circulation in large waterbodies is a management technique aimed to disrupt stratification of temperature and, consequently, 
to minimize the development of stagnant zones that may be subject to water quality problems (for instance, low levels of dissolved oxygen
or high concentrations of phytoplankton). 
For its operation, a set of flow pumps take water from the upper layers by means of collectors and inject it into the bottom layers, setting up 
a recirculation pattern that prevents stratification by means of a forced mixing of water. 
One of the main problems of the temperature stratification is related to algal blooms produced in the upper layers due to high temperature and solar radiation. 
However, if we circulate water from the bottom layers (where the temperature is lower) to the upper layers, we can mitigate this negative effect. 
Further details and remarks on several issues related to the optimal
design and control of water artificial circulation techniques have been analyzed by the authors in their recent work \cite{mcrf}.

Convective heat transfer has been the subject of an intensive mathematical research in last five decades
(ranging, for instance, from the pioneering works on the Boussinesq system of Joseph \cite{Jos} in the 1960s to the present).
Among the recent contributions we must mention, for instance, some papers devoted to study related problems 
in the steady case \cite{Benes,Kov}, the analysis a time-dependant case, but not including convective phenomena, \cite{Metzger1},
and some numerical approaches \cite{Bahl,Lamb}. Nevertheless, after an exhaustive search
we have not been able to find in the mathematical literature the analysis of the particular problem arising in the setting of our water recirculation model:
a coupled problem linking a heat equation with mixed nonlinear boundary conditions 
to a modified Navier-Stokes equation following the Smagorinsky model of turbulence.
Thus, the present work deals with the mathematical analysis and the numerical resolution of this 
heat transfer problem with specific boundary conditions related to water artificial circulation in a body of water (for instance, a lake or a reservoir). 
The main difficulties in the study of this problem lie in the nonlinear boundary condition related to the solar irradiation 
on the surface, the relations between the water temperature in the collectors and the injectors, 
and the coupling between water temperature and water velocity due to convective effects. 
We use the Smagorinsky model of turbulence instead of other 
approaches, like the celebrated $k-\epsilon$ system, due to the fact that the 
modified Navier-Stokes equations following the Smagorinsky model of turbulence 
present very interesting properties from a mathematical viewpoint, in particular 
the uniqueness of solution and its additional regularity.

The organization of this paper is as follows: First we introduce a 
well-posed formulation of the physical problem and present a rigorous 
definition of a solution for the problem. In the central part of the paper 
we prove the existence of this solution, and in the final part we propose 
a numerical algorithm for its resolution, showing several computational 
tests for a realistic example. At the end of the paper we include an 
appendix with several results for a general heat equation with 
an advective term and mixed boundary conditions of diffusive, convective 
and radiative type. These results, as far as we know, are new since we 
use techniques that will allow us to treat the low regularity of the time 
derivative of the solution. This lack of regularity represents, together 
with the nonlinearity of the boundary conditions, the main difficulty.

\section{Mathematical formulation of the problem}

\begin{figure}[ht]
\centering
\includegraphics[scale=1.]{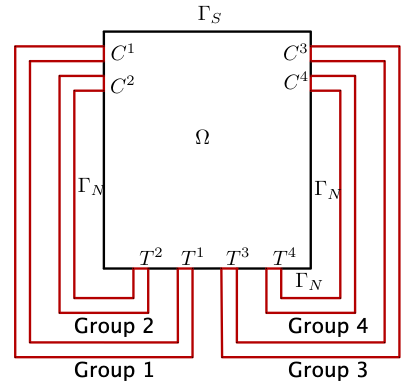}
\caption{Geometrical configuration of an example domain $\Omega$ with $N_{CT}=4$ collector/injector pairs,
showing the different boundary sections: $\Gamma_S$,  $\Gamma_C =\cup_{k=1}^{4} C^k$,  
$\Gamma_T=\cup_{k=1}^{4} T^k$, and $\Gamma_N$.}
\label{figure1}
\end{figure}

In this section we present in detail the three-dimensional mathematical model under study. 
So, we consider a convex domain $\Omega \subset \mathbb{R}^3$ 
(representing the waterbody) whose boundary surface $\partial \Omega$ can be split 
into four smooth enough, disjoint sections: $\Gamma_S$, $\Gamma_C$, 
$\Gamma_T$ and $\Gamma_N$, in such a way that $\partial \Omega=
\Gamma_S \cup \Gamma_C \cup \Gamma_T \cup \Gamma_N$. 
Subset $\Gamma_S$ represents the part of the boundary in contact with air, 
$\Gamma_C$ is the part of the boundary where the collectors 
are located, $\Gamma_T$ is the part of the boundary where the 
injectors are located, and $\Gamma_N$ stands for the rest of the boundary. 
We will suppose that each collector is linked to an injector by means 
of a pumped pipeline, and we also assume that there exist $N_{CT}$  
collector/injector pairs $\{ (C^k, T^k) \}_{k=1}^{N_{CT}}$. Therefore, $\Gamma_C=\cup_{k=1}^{N_{CT}} C^k$, and 
$\Gamma_T=\cup_{k=1}^{N_{CT}} T^k$. In Fig.~\ref{figure1} we 
can see a schematic geometrical configuration of a rectangular domain $\Omega$ 
for a particular case of $N_{CT}=4$ collector/injector pairs.

As above commented, we suppose the boundary $\partial \Omega$ regular enough to assure the existence of elements 
$\varphi^k,\, \widetilde{\varphi}^k \in H^{3/2}(\partial \Omega)$, 
$k=1,\ldots,N_{CT}$, satisfying the following assumptions (corresponding
to suitable regularizations of the indicator functions of $T^k$ and $C^k$, respectively):
\begin{itemize}
\item $\varphi^k(\mathbf{x}),\,\widetilde{\varphi}^k(\mathbf{x})\geq 0$, 
a.e. $\mathbf{x} \in \partial \Omega$, 
\item $\varphi^k(\mathbf{x})=0$, a.e. $\mathbf{x} \in \partial \Omega 
\setminus T^k$, and $\displaystyle \int_{T^k} \varphi^k(\mathbf{x}) \, d \gamma =\mu(T^k)$,
\item $\widetilde{\varphi}^k(\mathbf{x})=0$, a.e. $\mathbf{x} \in \partial \Omega 
\setminus C^k$, and $\displaystyle \int_{C^k} \widetilde{\varphi}^k(\mathbf{x})\, d \gamma = 
\mu(C^k)$,
\end{itemize}
where $\mu(S)$ represents the area measure of any set $S\subset \partial \Omega$.

We denote by $\theta(\mathbf{x},t)$ (measured in K) the solution of the following convection-diffusion partial 
differential equation with nonhomogeneous, nonlinear, mixed boundary conditions:
\begin{equation} \label{eq:system1}
\left\{\begin{array}{l}
\displaystyle 
\frac{\partial \theta}{\partial t}+ \mathbf{v} \cdot \nabla \theta - 
\nabla \cdot (K \nabla \theta) =0 \quad \mbox{in}\; \Omega \times ]0,T[, \\ 
\displaystyle
\theta =  \phi_{\theta} \quad \mbox{on} \; \Gamma_T\times ]0,T[, \\ 
\displaystyle
K \frac{\partial \theta}{\partial \mathbf{n}} = 0 \quad \mbox{on} \; \Gamma_C \times ]0,T[, \vspace{0.1cm}\\ 
\displaystyle
K \frac{\partial \theta}{\partial \mathbf{n}}= b_1^N (\theta_N-\theta) \quad \mbox{on} \; \Gamma_N \times ]0,T[, \vspace{0.1cm}\\ 
\displaystyle
K \frac{\partial \theta}{\partial \mathbf{n}} = b_1^S (\theta_S-\theta) +
b_2^S (T_r^4-|\theta|^3\theta) \quad \mbox{on}\; \Gamma_S \times ]0,T[,  \\ 
\theta(0)=\theta^0  \quad \mbox{in} \; \Omega,
\end{array}\right.
\end{equation}
where Dirichlet boundary condition $\phi_{\theta}$ is given by expression:
\begin{equation} \label{eq:defphi}
\phi_{\theta}(\mathbf{x},t)=\sum_{k=1}^{N_{CT}} \varphi^k(\mathbf{x})
\int_{-T}^T \rho_{\epsilon}(t-\epsilon-s) \, \gamma_{\theta}^k(s) \,  ds
\end{equation}
with, for each $k=1,\ldots,N_{CT}$, 
\begin{equation} \label{eq:gammatheta}
\gamma_{\theta}^k(s)=\left\{
\begin{array}{lcl}
\displaystyle \frac{1}{\mu(C^k)} \int_{C^k} \theta^0 \, d \gamma & \mbox{ if } & s \leq 0,  \\
\displaystyle \frac{1}{\mu(C^k)} \int_{C^k} \theta(s) \, d \gamma & \mbox{ if } & s > 0,
\displaystyle
\end{array}\right.
\end{equation}
representing the mean temperature of water in the collector $C_k$, and with the weight function 
$\rho_{\epsilon}$ defined by:
\begin{equation} \label{eq:rhoeps}
\rho_{\epsilon}(t)=\left\{\begin{array}{lcl}
\displaystyle \frac{c}{\epsilon} \, \exp\left(\frac{t^2}{t^2-\epsilon^2}\right) & \mbox{ if } &
|t| < \epsilon,  \\ 
\displaystyle 0 & \mbox{ if } & |t| \geq \epsilon,
\end{array}\right.
\end{equation}
for $c \in \mathbb{R}$ the positive constant satisfying the unitary condition:
\begin{equation}
\displaystyle \int_{\mathbb{R}} \rho_1(t) \, dt =1.
\end{equation}
In other words, we are assuming that the mean temperature of water at each injector $T_k$ is a 
weighted average in time of the mean temperatures of water at its corresponding collector $C_k$. 
In order to obtain the mean temperature at each injector, we convolute the 
mean temperature at the collector with a smooth function with support 
in $(t-2\epsilon,t)$. In this way, we have that the temperature in the injector 
only depends on the mean temperature in the collector in the time interval 
$(t-2\epsilon,t)$. Parameter $0<\epsilon<T$ represents, in a certain 
sense, the technical characteristics of the pipeline that define the stay time of water in the pipe. 
We also suppose that there is not heat transfer thought the walls of the pipelines (that is, they are isolated). 

Moreover, 
\begin{itemize}
\item $T>0  \ ({\rm s})$ is the length of the time interval.
\item $\mathbf{n}$ is the unit outward normal vector to the boundary $\partial \Omega$.
\item $K>0 \ ({\rm m}^2 \, {\rm s}^{-1})$ is the thermal diffusivity 
of the fluid: $K=\frac{\alpha}{\rho\, c_p}$, where 
$\alpha\ ({\rm W} \, {\rm m}^{-1} \, {\rm K}^{-1})$ is the thermal conductivity, 
$\rho\ ({\rm g} \, {\rm m}^{-3})$ is the density, and 
$c_p\ ({\rm W} \, {\rm s} \, {\rm g}^{-1} \, {\rm K}^{-1})$ is the specific heat 
capacity of water.
\item $b_1^K\geq 0 \ ({\rm m} \, {\rm s}^{-1})$, for $K=N,S$, are the coefficients 
related to convective heat transfer through the boundaries $\Gamma_N$ 
and $\Gamma_S$, obtained from the relation $\rho \, c_p\, b_1^K=h^K$, 
where $h^K\geq 0\ ({\rm W} \, {\rm m}^{-2} \, {\rm K}^{-1})$ 
are the convective heat transfer coefficients on each surface.  
These coefficients are relevant in the convective heat transfer flux through the frontiers 
$\Gamma_S$ and $\Gamma_N$, $b_1^K(\theta_K-\theta)$, $K=N,S$.
\item $b_2^S>0\ ({\rm m} \, {\rm s} \, {\rm K}^{-3})$ is the coefficient 
related to radiative heat transfer through the boundary $\Gamma_S$, given by
$b_2^S=\frac{\sigma_B\, \varepsilon}{\rho\, c_p}$, where 
$\sigma_B\ ({\rm W} \, {\rm m}^{-2} \, {\rm K}^{-4})$ 
is the Stefan-Boltzmann constant and $\varepsilon$ is the emissivity. 
This coefficient is fundamental in the radiative flux through the frontier 
$\Gamma_S$ (see, for instance, the classical reference \cite{Incropera} for 
a complete description of this type of boundary conditions).
\item $\theta^0\geq 0 \ ({\rm K})$ is the initial temperature. 
\item $\theta_S,\, \theta_N \geq 0 \ ({\rm K})$ are the temperatures 
related to convection heat transfer on the surfaces $\Gamma_S$ 
and $\Gamma_N$. 
\item $T_r \geq 0\ ({\rm K})$ is the radiation temperature on the 
surface $\Gamma_S$, derived from expression
$\sigma_B \, \varepsilon\, T_r^4 = (1-a) R_{sw,net}+R_{lw,down}$,
where $a$ is the albedo, $R_{sw,net}\ ({\rm W} \, {\rm m}^{-2})$ 
denotes the net incident shortwave radiation on the surface $\Gamma_R$,
and $R_{lw,down} \ ({\rm W} \, {\rm m}^{-2})$ denotes 
the downwelling longwave radiation.
\end{itemize}

Finally, $\mathbf{v}(\mathbf{x},t) \ ({\rm m} \, {\rm s}^{-1})$ is the water velocity, solution 
of a modified Navier-Stokes equations following a Smagorinsky model of turbulence:
\begin{equation} \label{eq:system2}
\left\{\begin{array}{l}
\displaystyle
\frac{\partial \mathbf{v}}{\partial t} + \nabla \mathbf{v} \, \mathbf{v} - \nabla \cdot \Xi(\mathbf{v})
+\nabla p = \alpha_0 (\theta-\theta^0) \mathbf{a}_g \quad \mbox{in} \; \Omega \times ]0,T[,  \\
\displaystyle
\nabla \cdot \mathbf{v}=0 \quad \mbox{in} \; \Omega \times ]0,T[, \\
\displaystyle
\mathbf{v}=\boldsymbol{\phi}_{\mathbf{g}} \quad  \mbox{on}\; \partial \Omega\times ]0,T[,  \\
\displaystyle
\mathbf{v}(0)=\mathbf{v}^0 \quad \mbox{in}\; \Omega,
\end{array}\right.
\end{equation}
where $\mathbf{a}_g=g\, \mathbf{e}_3$ $({\rm m} \, {\rm s}^{-1})$ is the 
gravity acceleration,  
$\alpha_0=-\frac{1}{\rho} \frac{\partial \rho}{\partial \theta}$ $({\rm K}^{-1})$ is the 
thermic expansion coefficient. We must remark here that we are assuming 
the thermodynamic process to be close to an initial equilibrium state that we denote with the zero subscript, so $\rho=\rho_0-\alpha_0 \rho_0 (\theta-\theta_0)+ 
k_0 \rho_0(p-p_0)+o(\theta-\theta_0)+o(p-p_0)$, where $k_0=\frac{1}{\rho} 
\frac{\partial \rho}{\partial p}$ is the isothermal compressibility 
coefficient. The details of this approach (known as the Boussinesq model for 
natural convection) can be consulted, for instance, in Section 10.7 of 
\cite{bermudez}. Finally, $\mathbf{v}^0$ is the initial 
velocity, and boundary field $\boldsymbol{\phi}_{\mathbf{g}}$ is the element given by:
\begin{equation} \label{eq:phig}
\boldsymbol{\phi}_{\mathbf{g}}(t,\mathbf{x})
=\sum_{k=1}^{N_{CT}}g^k(t) \left[\frac{\varphi^k(\mathbf{x})}{\mu(T^k)}
 -\frac{\widetilde{\varphi}^k(\mathbf{x})}{\mu(C^k)} 
\right] \mathbf{n}.
\end{equation}
with, for each $k=1,\ldots,N_{CT}$, 
$g^k(t) \in H^1(0,T)$, representing the volumetric flow rate by 
pump $k$ at each time $t$ ($g^k(t)>0$, $\forall t \in ]0,T[$, and $g^k(0)=0$). 
The turbulence term $\Xi(\mathbf{v})$ is given by:
\begin{equation}
\Xi(\mathbf{v})=\left. \frac{\partial D(e)}{\partial e} \right|_{e=e(\mathbf{v})},
\ \mbox{ with } e(\mathbf{v})=\frac{1}{2} \left( \nabla \mathbf{v} + \nabla \mathbf{v}^t\right),
\end{equation}
where $D$ is a potential function (for instance, in the particular case 
of the classical Navier-Stokes equations,
$D(e)=\nu \left[e:e\right]$, with $\nu$ $({\rm m}^2 \, {\rm s}^{-1})$ the kinematic viscosity of the water, and, consequently,
$\Xi(\mathbf{v}) = 2 \nu\, e(\mathbf{v})$).
However, in our case, the Smagorinsky model, the potential function is defined as \cite{lady1}:
\begin{equation} \label{eq:doperador}
D(e)=\nu \left[e:e\right] + \frac{2}{3} \nu_{tur} \left[e:e\right]^{3/2},
\end{equation}
where $\nu_{tur}$ $({\rm m}^2)$ is the turbulent viscosity.  So, for the Smagorinsky case,
\begin{equation}
\begin{array}{rcl}
\displaystyle
\Xi(\mathbf{v}) &=&
\displaystyle \left. \frac{\partial D(e)}{\partial e} \right|_{e=\epsilon(\mathbf{v})} = \displaystyle
2 \nu \, \epsilon(\mathbf{v}) + 2 \nu_{tur}  \left[\epsilon(\mathbf{v}):\epsilon(\mathbf{v})\right]^{1/2}
\epsilon(\mathbf{v}) \\
&=& \displaystyle
\left(2 \nu+ 2 \nu_{tur}  \left[\epsilon(\mathbf{v}):\epsilon(\mathbf{v})\right]^{1/2} \right) \epsilon(\mathbf{v}) = \displaystyle
\beta(\epsilon(\mathbf{v})) \, \epsilon(\mathbf{v}),
\end{array}
\end{equation}
with $\beta(\epsilon(\mathbf{v}))=
2 \nu+ 2 \nu_{tur}  \left[\epsilon(\mathbf{v}):\epsilon(\mathbf{v})\right]^{1/2} $.

System (\ref{eq:system2}) has been recently studied by the authors in \cite{fran7}, where it is demonstrated the existence and the uniqueness of solution
for a recirculation model based in the modified Navier-Stokes equations. In 
the present work we will use some of the results shown in \cite{fran7} in order to 
prove the existence of solution for the coupled problem (\ref{eq:system1}) and (\ref{eq:system2}). The main difficulties of this work relies in the coupling of a heat 
equation with nonlinear boundary value terms and the modified Navier-Stokes system. 
The nonlinear terms do not allow us to work with regular solutions for the heat equation, which forces us to use more sophisticated techniques in order to demonstrate existence and uniqueness of solution.

\section{The concept of solution}

We start this section defining the functional spaces used in the definition of solution 
for the system (\ref{eq:system1}) and (\ref{eq:system2}). So, for the water 
temperature we consider:
\begin{equation}
\begin{array}{rcl}
X_1&=&\{\theta \in H^1(\Omega):\; \theta_{|_{\Gamma_S}}\in L^5(\Gamma_S)\}, \\ 
\widetilde{X}_1 &=& \{\theta \in X_1:\; \theta_{|_{\Gamma_T}}=0\},
\end{array}
\end{equation}
and we define the following norm associated to above space $X_1$:
\begin{equation}
\|\theta\|_{X_1}=\|\theta\|_{H^1(\Omega)}+\|\theta\|_{L^5(\Gamma_S)}.
\end{equation}
We have that $X_1$ is a reflexive separable Banach space 
(cf. Lemma 3.1 of \cite{Zolesio1}) and $\widetilde{X}_1 \subset L^2(\Omega) 
\subset \widetilde{X}_1'$ is an evolution triple. For the water velocity we 
consider:
\begin{equation}
\begin{array}{rcl}
\displaystyle
\mathbf{X}_2&=&
\displaystyle
\left\{\mathbf{v}\in [W^{1,3}(\Omega)]^3:\;
\nabla \cdot \mathbf{v}=0 , \;  \mathbf{v}_{|_{\Gamma_S \cup \Gamma_N}}= \boldsymbol{0}
\right\} , \\
\displaystyle
\widetilde{\mathbf{X}}_2&=& \displaystyle \left\{\mathbf{v}\in [W^{1,3}(\Omega)]^3:\;
\nabla \cdot \mathbf{v}=0, \; \mathbf{v}_{|_{\partial \Omega}}= \boldsymbol{0}
\right\} .
\end{array}
\end{equation}

In order to define an appropriate space for the solution of problems (\ref{eq:system1}) and (\ref{eq:system2}), 
we consider, for a Banach space $V_1$ and a locally convex space $V_2$ such that $V_1\subset V_2$, 
the following Sobolev-Bochner space (cf. Chapter 7 of \cite{Roubicek1}), for $1\leq p,q \leq \infty$:
\begin{equation}
W^{1,p,q}(0,T;V_1,V_2)=\left\{u \in L^p(0,T;V_1):\; \frac{du}{dt} \in L^q(0,T;V_2)\right\},
\end{equation}
where $\frac{du}{dt} $ denotes the derivative of $u$ in the sense of distributions. It is well known that, if both $V_1$ and $V_2$ are Banach 
spaces, then $W^{1,p,q}(0,T;V_1,V_2)$ is also a Banach space endowed 
with the norm $\|u\|_{W^{1,p,q}(0,T;V_1,V_2)}=\|u\|_{L^p(0,T;V_1)}+
\left\| \frac{du}{dt}  
\right\|_{L^{q}(0,T;V_2)}$.

Then, we define the following spaces that will be used in 
the mathematical analysis of system (\ref{eq:system1}):
\begin{equation} \hspace*{-1.5cm}
\begin{array}{rcl}
\displaystyle
W_1&=&\displaystyle
\{\theta \in W^{1,2,5/4}(0,T;X_1,X_1'):\\ && \displaystyle \theta_{|_{\Gamma_S}} \in 
L^5(0,T;L^5(\Gamma_S))\}\cap L^{\infty}(0,T;L^2(\Omega)), \\
\widetilde{W}_1 &=& \displaystyle
\{\theta \in W^{1,2,5/4}(0,T;\widetilde{X}_1,\widetilde{X}_1'):\\ && 
\displaystyle
\theta_{|_{\Gamma_S}} \in 
L^5(0,T;L^5(\Gamma_S))\}\cap L^{\infty}(0,T;L^2(\Omega)),
\end{array}
\end{equation}
and, for the system (\ref{eq:system2}), we define 
\begin{equation}
\begin{array}{rcl}
\displaystyle
\mathbf{W}_2&=&
\displaystyle
W^{1,\infty,2}(0,T;\mathbf{X}_2,[L^2(\Omega)]^3) \cap 
\mathcal{C}([0,T];\mathbf{X}_2), \\
\displaystyle
\widetilde{\mathbf{W}}_2&=&
\displaystyle
W^{1,\infty,2}(0,T;\widetilde{\mathbf{X}}_2,[L^2(\Omega)]^3) \cap 
\mathcal{C}([0,T];\widetilde{\mathbf{X}}_2).
\end{array}
\end{equation}

\begin{hypo} \label{hypo1} We will assume the following hypotheses for coefficients and data of the problem:
\begin{enumerate}
\item[(a)] $\theta^0 \in X_2$
\item[(b)] $\theta_S \in L^2(0,T;L^2(\Gamma_S))$
\item[(c)] $\theta_N \in L^2(0,T;L^2(\Gamma_N))$
\item[(d)] $T_r \in L^5(0,T;L^5(\Gamma_S))$
\item[(e)]  $\mathbf{v}^0 \in \left[H_{\sigma}^2(\Omega)\right]^3=\{\mathbf{v} \in 
[H^2(\Omega)]^3 :\; \nabla \cdot \mathbf{v}=0, \ \mathbf{v}_{|_{\partial \Omega}}=\boldsymbol{0}\}\subset 
\widetilde{\mathbf{X}}_2$
\item[(f)] $g^k \in H^1(0,T)$ with $g^k(0)=0$, $\forall k=1,\ldots,N_{CT}$
\end{enumerate}
\end{hypo}

\begin{rem} In order to define in a rigorous way the concept of solution, we will need to extend Dirichlet conditions
of $\theta$ and $\mathbf{v}$ to the whole domain $\Omega$. 

So, for water velocity $\mathbf{v}$, thanks to Lemma 2 of \cite{fran7}, for 
each $\mathbf{g} \in [H^1(0,T)]^{N_{CT}}$, there exists an element $\boldsymbol{\zeta}_{\mathbf{g}} \in 
W^{1,2,2}(0,T;\left[H_{\sigma}^2(\Omega)\right]^3,
\left[H_{\sigma}^2(\Omega)\right]^3)$ such that ${\boldsymbol{\zeta}_{\mathbf{g}}}_{|_{\partial \Omega}}=
\boldsymbol{\phi}_{\mathbf{g}}$, with $\boldsymbol{\phi}_{\mathbf{g}}$ defined by (\ref{eq:phig}). Besides, 
by Lemma 3 of \cite{fran7}, $W^{1,2,2}(0,T;\left[H_{\sigma}^2(\Omega)\right]^3,$
$\left[H_{\sigma}^2(\Omega)\right]^3) \subset W^{1,\infty,2}(0,T;[H^2(\Omega)]^3,[H^2(\Omega)]^3) 
\cap \mathcal{C}([0,T];[H^2(\Omega)]^3)$ and, then, 
we can use this element to reformulate the original problem for $\mathbf{v}$ as an homogeneous Dirichlet boundary condition one. 

For water temperature $\theta$ we can proceed in an analogous way and prove that 
there exists an extension that allows us to reformulate the problem for $\theta$ as one with homogeneous boundary conditions. 
\end{rem}

\begin{lem} \label{lemma2} We have that the following operator is compact
\begin{equation}
\begin{array}{rcl}
R_{\mathbf{h}}:[L^2(0,T)]^{N_{CT}} & \rightarrow & 
W^{1,2,2}(0,T;H^{2}(\Omega),H^2(\Omega))\\ 
\mathbf{h} & \rightarrow & R_{\mathbf{h}}(\mathbf{h} )=\zeta_{\mathbf{h} },
\end{array}
\end{equation}
where:
\begin{equation}
\zeta_{\mathbf{h} }(\mathbf{x},t)=\sum_{k=1}^{N_{CT}} \beta_0(\varphi^k(\mathbf{x}) )
\int_{-T}^T \rho_{\epsilon} (t-\epsilon-s) 
\gamma^k_{\mathbf{h} }(s) \, ds,
\end{equation}
with $\gamma^k_{\mathbf{h}}(s) \in L^2(-T,T)$, for $k=1,2,\ldots,N_{CT}$, 
defined by:
\begin{equation}
\gamma^k_{\mathbf{h}}(s)=\left\{
\begin{array}{lcl}
\displaystyle \frac{1}{\mu(C^k)} \int_{C^k} \theta^0 \, d \gamma & \mbox{ if } & s \leq 0,   \\
\displaystyle{h^k}(s) & \mbox{ if } & s > 0,
\end{array}\right.
\end{equation}
and $\beta_0: u \in H^{3/2}(\partial \Omega) \rightarrow \beta_0(u) \in H^2(\Omega)$ the right inverse of the classical trace operator $\gamma_0$
(that is, $(\gamma_0 \circ \beta_0)(u) =u$.) 

We also have that there exists 
a constant $C_1$, that depends continuously on the space-time configuration of 
our computational domain and $\theta^0$, such that:
\begin{equation} \label{eq:acotapf0}
\|\zeta_{\mathbf{h}}\|_{W^{1,2,2}(0,T;H^2(\Omega),H^2(\Omega))} \leq C_1(\theta^0)(1+
\|\mathbf{h}\|_{[L^2(0,T)]^{N_{CT}}}).
\end{equation}
\end{lem}

\begin{pf} Let $\{\mathbf{h}_n\}_{n \in \mathbb{N}}$ be a bounded sequence in 
$[L^2(0,T)]^{N_{CT}}$. Then, taking subsequences if necessary, we have that 
$\mathbf{h}_n \rightharpoonup \mathbf{h}$ weakly in $[L^2(0,T)]^{N_{CT}}$. 
We also have that $\gamma_{\mathbf{h}_n} \rightharpoonup \gamma_{\mathbf{h}}$ 
weakly in $[L^2(-T,T)]^{N_{CT}}$ and, if we denote by 
\begin{equation}
b_{\mathbf{h}}^k(t)=\int_{-T}^T \rho_{\epsilon} (t-\epsilon-s) 
\gamma^k_{\mathbf{h}}(s) \, ds,
\end{equation}
we obtain that $b_{\mathbf{h}_n}^k(t) \rightarrow b_{\mathbf{h}}^k(t)$ pointwise a.e. $t \in [0,T]$, 
$\forall k=1,\ldots,N_{CT}$. Thus, the sequence $\{\mathbf{b}_{\mathbf{h}_n}\}_{n \in \mathbb{N}}$ is  
bounded by a function in $[L^2(0,T)]^{N_{CT}}$, so we have the strong convergence in 
$[L^2(0,T)]^{N_{CT}}$. We can repeat the same argument with the time derivative of 
$\{\mathbf{b}_{\mathbf{h}_n}\}_{n \in \mathbb{N}}$, obtaining that
$\mathbf{b}_{\mathbf{h}_n} \rightarrow \mathbf{b}_{\mathbf{h}}$ strongly in $H^1(0,T)$, thus 
$\zeta_{\mathbf{h}_n} \rightarrow \zeta_{\mathbf{h}}$ strongly in $W^{1,2,2}(0,T;H^2(\Omega), 
H^2(\Omega))$. Finally, by the properties of the operator $\beta_0$ we have that
\begin{equation}
\|\beta_0(\varphi^k )\|_{H^{2}(\Omega)} \leq C \|\varphi^k\|_{H^{3/2}(\partial \Omega)}, \;
k=1,\ldots,N_{CT}.
\end{equation}
So, thanks to the regularity of function $\rho_{\epsilon}$, it is clear that $b_{\mathbf{h}}^k \in H^1(0,T)$, 
and
\begin{equation}
\|b_{\mathbf{h}}^k\|_{H^1(0,T)} \leq C \|\gamma_{\mathbf{h}}^k\|_{L^2(-T,T)}.
\end{equation}
In the other hand, 
\begin{equation}
\begin{array}{rcl}
\displaystyle
\|\gamma_{\mathbf{h}}^k\|_{L^2(-T,T)}^2&=&
\displaystyle
\int_{-T}^T \left(\gamma_{\mathbf{h}}^k(s) \right)^2\, ds \\ &=&
\displaystyle
\frac{1}{\mu(C^k)^2} \int_{-T}^0 \left(\int_{C^k} \theta^0 \, d \gamma  \right)^2+
 \int_0^T \left(  h^k(s) \right)^2\, ds
\\ &\leq &
\displaystyle 
\frac{ T}{\mu(C^k)} \|\theta^0\|_{L^2(\partial \Omega)}^2+
\|h^k\|^2_{L^2(0,T)}.
\end{array}
\end{equation}
Thereby, we have the following inequality:
\begin{eqnarray*}
\displaystyle
\|\zeta_{\mathbf{h}}\|_{W^{1,2,2}(0,T;H^2(\Omega),H^2(\Omega))} & \leq & C\left( 
\|\theta^0\|_{L^2(\partial \Omega)}+\|\mathbf{h}\|_{[L^2(0,T)]^{N_{CT}}}\right)
\sum_{k=1}^{N_{CT}} \mu(T^k) \\
& \leq & C_1(\theta^0) \left( 1+\|\mathbf{h}\|_{[L^2(0,T)]^{N_{CT}}} \right),
\end{eqnarray*}
where $C(\theta^0)$ is a positive constant than depends continuously on 
the spatial-time configuration of our computational domain and on the initial temperature. 
Moreover, it is worthwhile remarking here that we can make this constant as small as we want 
by considering data appropriately. \hfill $\blacksquare$
\end{pf}

Finally, the following technical lemma will be necessary in order to guaranty that the sum of an element of 
$W_1$ plus an element of $W^{1,2,2}(0,T;H^{2}(\Omega),H^{2}(\Omega))$ makes sense. 

\begin{lem} \label{lemma3} We have that the following inclusion is compact:
\begin{equation}
W^{1,2,2}(0,T;H^{2}(\Omega),H^{2}(\Omega)) \subset \subset W_1. 
\end{equation}
\end{lem}

Now, we define the concept of solution for coupled system (\ref{eq:system1}) and 
(\ref{eq:system2}) in terms of homogeneous Dirichlet systems, so we 
establish the following notations:
\begin{itemize}
\item $\xi=\theta - \zeta_{\mathbf{h}}\in \widetilde{W}_1$, 
with $\displaystyle \zeta_{\mathbf{h}} \in W^{1,2,2}(0,T;H^{2}(\Omega),H^{2}(\Omega))$ 
the extension obtained from Lemma \ref{lemma2}, where:
\begin{equation}
h^k(s)=\frac{1}{\mu(C^k)} \int_{C^k} \theta(s) \, d \gamma, \quad k=1,2,\ldots,N_{CT}.
\end{equation}
\item $\mathbf{z}=\mathbf{v} - \boldsymbol{\zeta}_{\mathbf{g}} 
\in \widetilde{\mathbf{W}}_2$, 
with $\boldsymbol{\zeta}_{\mathbf{g}} \in W^{1,2,2}(0,T;[H_{\sigma}^2(\Omega)]^3,[H_{\sigma}^2(\Omega)]^3)$ 
the extension of the trace given in Lemma 2 of \cite{fran7}.
\end{itemize}
Thus, using above notations, we can reformulate the state system
 (\ref{eq:system1}) and (\ref{eq:system2}) in the following way:

\begin{equation}\label{eq:system1b}
\left\{\begin{array}{l}
\displaystyle \frac{\partial \xi}{\partial t}+\mathbf{v} \cdot \nabla \xi - \nabla \cdot (K \nabla \xi)  \\ 
\displaystyle \quad = -\frac{\partial \zeta_{\mathbf{h}}}{\partial t}-\mathbf{v} \cdot \nabla \zeta_{\mathbf{h}} + 
\nabla \cdot (K \nabla \zeta_{\mathbf{h}}) \quad \mbox{in} \; \Omega \times (0,T), \\ 
\displaystyle 
\xi=0 \quad \mbox{on} \; T^k \times (0,T),  \quad  \mbox{for} \; k=1,\ldots,N_{CT}, \\ 
\displaystyle K \frac{\partial \xi}{\partial \mathbf{n}}= -K  \frac{\partial \zeta_{\mathbf{h}}}{\partial \mathbf{n}} 
\quad \mbox{on} \; C^k \times (0,T), \quad \mbox{for} \; k=1,\ldots,N_{CT}, \\ 
\displaystyle K \frac{\partial \xi}{\partial \mathbf{n}}=b_1^N\big(
\theta_N-\zeta_{\mathbf{h}}-\frac{K}{b_1^N} \frac{\partial \zeta_{\mathbf{h}}}{\partial \mathbf{n}}-\xi
\big) \quad \mbox{on}\; \Gamma_N \times (0,T), \\ 
\displaystyle K \frac{\partial \xi}{\partial \mathbf{n}}=b_1^S\big(
\theta_S-\zeta_{\mathbf{h}}-\frac{K}{b_1^S} \frac{\partial \zeta_{\mathbf{h}}}{\partial \mathbf{n}}-\xi \big) \\
\displaystyle \quad +b_2^S \big(T_r^4-|\xi+\zeta_{\mathbf{h}}|^3(\xi+\zeta_{\mathbf{h}})\big) 
\quad \mbox{on}\; \Gamma_S \times (0,T), \\ 
\displaystyle \xi(0)=\theta^0-\zeta_{\mathbf{h}}(0) \quad \mbox{in} \; \Omega.
\end{array}\right.
\end{equation}
\begin{equation} \label{eq:system2b}
\left\{
\begin{array}{l}
\displaystyle \frac{\partial \mathbf{z}}{\partial t} + \nabla (\boldsymbol{\zeta}_{\mathbf{g}} +
\mathbf{z}) \mathbf{z} +\nabla \mathbf{z} \boldsymbol{\zeta}_{\mathbf{g}} \\ 
\displaystyle  \quad -\text{div}  \left( 2 \nu \epsilon(\mathbf{z})+
2 \nu_{tur} \int_{\Omega} \left[\epsilon(\boldsymbol{\zeta}_{\mathbf{g}}+\mathbf{z}):
\epsilon(\boldsymbol{\zeta}_{\mathbf{g}}+\mathbf{z}) \right]^{1/2} \epsilon(\boldsymbol{\zeta}_{\mathbf{g}}+\mathbf{z}) \right) \\ 
\displaystyle \quad +\nabla p = \alpha_0(\theta-\theta^0)\, \mathbf{a}_{g}- \frac{\partial \boldsymbol{\zeta}_{\mathbf{g}}}{\partial t}
-\nabla \boldsymbol{\zeta}_{\mathbf{g}} \boldsymbol{\zeta}_{\mathbf{g}} 
\\ \quad \displaystyle
+2 \nu \nabla \cdot 
\epsilon(\boldsymbol{\zeta}_{\mathbf{g}})\quad \mbox{in} \; \Omega \times (0,T), \\
\displaystyle \mathbf{z}  = \mathbf{0} \quad \mbox{on} \; \partial \Omega \times (0,T),  \\
\displaystyle \mathbf{z}(0)=\mathbf{v}^0 \quad \mbox{in}\; \Omega.
\end{array}
\right.
\end{equation}
It is worthwhile remarking here that above system shows homogeneous Dirichlet 
boundary conditions and, consequently, we will be able to define the concept 
of solution of the original state systems (\ref{eq:system1}) and (\ref{eq:system2}) 
in terms of the modified state systems (\ref{eq:system1b}) and 
(\ref{eq:system2b}). It should be also noted that, in the case of equation 
(\ref{eq:system1}), the coupling terms in the Dirichlet boundary conditions are now 
transferred to the partial differential equation in system (\ref{eq:system1b}).

\begin{defn} \label{defsol} A pair $(\theta,\mathbf{v})\in W_1 \times \mathbf{W}_2$ 
is said to be a solution of problem (\ref{eq:system1}) and (\ref{eq:system2}) if
there exist elements $(\xi,\mathbf{z})\in \widetilde{W}_1\times \widetilde{\mathbf{W}}_2$ 
such that:
\begin{enumerate}
\item $\mathbf{v}=\boldsymbol{\zeta}_{\mathbf{g}} + \mathbf{z}$, with
$\displaystyle \boldsymbol{\zeta}_{\mathbf{g}} \in 
W^{1,2,2}(0,T;[H^2(\Omega)]^3,[H^2(\Omega)]^3)$
the reconstruction of the trace given in Lemma 2 of \cite{fran7}, and  
$\theta=\zeta_{\mathbf{h}}+\xi$, with $\displaystyle \zeta_{\mathbf{h}}\in 
W^{1,2,2}(0,T;H^{2}(\Omega),H^{2}(\Omega))$ the extension obtained 
in Lemma \ref{lemma2}, where:
\begin{equation}
h^k(s)=\frac{1}{\mu(C^k)} \int_{C^k} \theta(s) \, d \gamma, \quad k=1,2,\ldots,N_{CT}.
\end{equation}
\item $\mathbf{z}(0)=\mathbf{v}^0 $, and $ \xi(0)=\theta^0-\zeta_{\mathbf{h}}(0)$, a.e. 
$\mathbf{x} \in \Omega$.
\item $(\xi,\mathbf{z})$ verifies the following variational formulation:
\begin{equation} \label{eq:system3}
\hspace{-1.6cm} \begin{array}{r}
\displaystyle \int_{\Omega} \frac{\partial \xi}{\partial t} \eta \, d \mathbf{x} +
\int_{\Omega} \mathbf{v} \cdot \nabla \xi \eta \, d \mathbf{x} + 
K \int_{\Omega} \nabla \xi \cdot \nabla \eta \, d \mathbf{x} +  
b_1^N \int_{\Gamma_N} \xi \eta \, d \gamma \\ 
\displaystyle 
+b_1^S \int_{\Gamma_S} \xi \eta \, d \gamma  
\displaystyle 
+b_2^S \int_{\Gamma_S} |\xi+\zeta_{\mathbf{h}}|^3(\xi+\zeta_{\mathbf{h}}) \eta 
\, d \gamma = \int_{\Omega} H_{\mathbf{h}}  \eta \, d \mathbf{x} 
\\ 
\displaystyle 
+\int_{\Gamma_C} g_{\mathbf{h}}^{C} \eta \, d \gamma 
+b_1^N \int_{\Gamma_N}   g_{\mathbf{h}}^{N} \eta \, d \gamma 
\displaystyle
+b_1^S \int_{\Gamma_S}   g_{\mathbf{h}}^{S} \eta \, d \gamma
+b_2^S \int_{\Gamma_S} T_r^4 \eta \, d \gamma, 
\\ 
\displaystyle  \mbox{a.e.} \ t \in ]0,T[, \quad 
\forall \eta \in \widetilde{X}_1.
\end{array}
\end{equation}
\begin{equation} \label{eq:system4}
\hspace{-1.6cm} \begin{array}{r}
\displaystyle
\int_{\Omega} \frac{\partial \mathbf{z}}{\partial t} \cdot \boldsymbol{\eta} \,d\mathbf{x} +
\int_{\Omega} \nabla (\boldsymbol{\zeta}_{\mathbf{g}}+
\mathbf{z}) \mathbf{z} \cdot \boldsymbol{\eta} \, d \mathbf{x} +
\int_{\Omega} \nabla \mathbf{z} \boldsymbol{\zeta}_{\mathbf{g}} 
\cdot \boldsymbol{\eta} \, d \mathbf{x} 
\\ 
\displaystyle 
+ 2 \nu \int_{\Omega} \epsilon(\mathbf{z}) : \epsilon(\boldsymbol{\eta}) \, 
d \mathbf{x} \\
+ \displaystyle
2 \nu_{tur} \int_{\Omega} \left[\epsilon(\boldsymbol{\zeta}_{\mathbf{g}}+\mathbf{z}):
\epsilon(\boldsymbol{\zeta}_{\mathbf{g}}+\mathbf{z}) \right]^{1/2} \epsilon(\boldsymbol{\zeta}_{\mathbf{g}}+\mathbf{z}):
\epsilon(\boldsymbol{\eta}) \, d \mathbf{x} \\
\displaystyle
=\int_{\Omega} \mathbf{H}_{\mathbf{g}} \cdot \boldsymbol{\eta} \, d \mathbf{x},
\quad \mbox{a.e.} \ t \in ]0,T[, \quad \forall
 \boldsymbol{\eta} \in \widetilde{\mathbf{X}}_2,
\end{array}
\end{equation}
where
\begin{equation}
\begin{array}{rcl}
H_{\mathbf{h}}&=&\displaystyle 
\frac{\partial \zeta_{\mathbf{h}}}{\partial t}-\mathbf{v} \cdot \nabla \zeta_{\mathbf{h}} + 
\nabla \cdot (K \nabla \zeta_{\mathbf{h}}) \in L^2(0,T;L^2(\Omega)), \\
\displaystyle
g^C_{\mathbf{h}} &=&\displaystyle 
-K  \frac{\partial \zeta_{\mathbf{h}}}{\partial \mathbf{n}} \in L^2(0,T;L^2(\Gamma_C)) ,
\\ 
\displaystyle
g^N_{\mathbf{h}}&=&\displaystyle
\theta_N-\zeta_{\mathbf{h}}-\frac{K}{b_1^N} \frac{\partial \zeta_{\mathbf{h}}}{\partial \mathbf{n}}
\in L^2(0,T;L^2(\Gamma_N)), 
\\ 
\displaystyle
g^S_{\mathbf{h}}&=&\displaystyle
\theta_S-\zeta_{\mathbf{h}}-\frac{K}{b_1^S} \frac{\partial \zeta_{\mathbf{h}}}{\partial \mathbf{n}}
\in L^2(0,T;L^2(\Gamma_S)),
\\ 
\displaystyle
 \mathbf{H}_{\mathbf{g}}&=&\displaystyle
 \alpha_0 (\theta-\theta^0) \mathbf{a}_g-
\frac{\partial \boldsymbol{\zeta}_{\mathbf{g}}}{\partial t}
-\nabla \boldsymbol{\zeta}_{\mathbf{g}} \boldsymbol{\zeta}_{\mathbf{g}} \in L^2(0,T;[L^2(\Omega)]^3).
\end{array}
\end{equation}
\end{enumerate}
\end{defn}

\section{Existence of solution}

We will prove now that, under certain hypotheses over coefficients 
and data, there exists a unique solution for the system (\ref{eq:system1}) 
and (\ref{eq:system2}) in the sense of Definition \ref{defsol}. 
The procedure used here for demonstrating the existence of solution is based in the Schauder 
fixed point Theorem (cf. section 9.5 of \cite{Conway}) and 
is similar to one employed by the authors, for instance, in \cite{fran6}. 
The main difficulties in the present case lie in the coupling of the 
Dirichlet conditions for the water temperature, and in the nonlinear 
radiation terms. To overcome these difficulties we will need to define a 
compact extension for the nonhomogeneous Dirichlet conditions and 
to prove novel results for the heat equation with radiation boundary conditions (see Appendix \ref{app}).

So, we consider the following operator:
\begin{equation} \hspace{-.15cm} \label{puntofijo1}
\begin{array}{rcl}
\mathbf{M} : L^2(0,T;L^2(\Omega)) \times [L^2(0,T)]^{N_{CT}} & \rightarrow & 
L^2(0,T;L^2(\Omega)) \times [L^2(0,T)]^{N_{CT}} \\
(\theta^*,\mathbf{h}^*) & \rightarrow & \mathbf{M} (\theta^*,\mathbf{h}^*)= (\theta,\mathbf{h}) 
\end{array}
\end{equation}
where:
\begin{itemize}
\item $\mathbf{v} \in \mathbf{W}_2$ is such that $\mathbf{z}=\mathbf{v}-\boldsymbol{\zeta}_{\mathbf{g}}
\in \widetilde{\mathbf{W}}_2$ is the solution of the following problem:
\begin{equation} \label{eq:system6} \hspace{-.3cm}
\left\{
\begin{array}{l}
\displaystyle \frac{\partial \mathbf{z}}{\partial t} + \nabla (\boldsymbol{\zeta}_{\mathbf{g}} +
\mathbf{z}) \mathbf{z} +\nabla \mathbf{z} \boldsymbol{\zeta}_{\mathbf{g}}  \\ 
\quad \displaystyle  -\nabla \cdot \left(2 \nu \epsilon(\mathbf{z})+
2 \nu_{tur}\left[\epsilon(\boldsymbol{\zeta}_{\mathbf{g}}+\mathbf{z}):
\epsilon(\boldsymbol{\zeta}_{\mathbf{g}}+\mathbf{z}) \right]^{1/2} \epsilon(\boldsymbol{\zeta}_{\mathbf{g}}+\mathbf{z}) 
\right)  \\ 
\quad \displaystyle
+\nabla p = \alpha_0(\theta^*-\theta^0) \mathbf{a}_{g}-
\frac{\partial \boldsymbol{\zeta}_{\mathbf{g}}}{\partial t}
-\nabla \boldsymbol{\zeta}_{\mathbf{g}} \boldsymbol{\zeta}_{\mathbf{g}} 
+2 \nu \nabla \cdot 
\epsilon(\boldsymbol{\zeta}_{\mathbf{g}}) \ \mbox{in} \; \Omega \times ]0,T[, \\ 
\displaystyle
\mathbf{z}  = \mathbf{0} \ \mbox{on} \; 
\Gamma \times ]0,T[,  \\
\displaystyle
\mathbf{z}(0)=\mathbf{v}^0-\boldsymbol{\zeta}_{\mathbf{g}}(0) \ \mbox{in}\; \Omega,
\end{array}
\right.
\end{equation}
with $\boldsymbol{\zeta}_{\mathbf{g}}\in 
W^{1,2,2}(0,T;[H_{\sigma}^{2}(\Omega)]^3,[H_{\sigma}^{2}(\Omega)]^3)$ 
the extension obtained from the Lemma 2 of \cite{fran7}.
\item $\theta \in W_1$ is such that $\xi=\theta-\zeta_{\mathbf{h}^*} \in \widetilde{W}_1$ is the solution of:
\begin{equation} \label{eq:system5}
\left\{\begin{array}{l}
\displaystyle
\frac{\partial \xi}{\partial t}+\mathbf{v} \cdot \nabla \xi - \nabla \cdot (K \nabla \xi) = 
-\frac{\partial \zeta_{\mathbf{h}^*}}{\partial t}-\mathbf{v} \cdot \nabla \zeta_{\mathbf{h}^*} \\
\qquad + \nabla \cdot (K \nabla \zeta_{\mathbf{h}^*}) \quad \mbox{in} \; \Omega \times ]0,T[, 
\\ 
\displaystyle 
\xi=0 \quad \mbox{on}\; T^k \times ]0,T[, \ \mbox{for}\, k=1,\ldots,N_{CT},
\\ 
\displaystyle
K \frac{\partial \xi}{\partial \mathbf{n}}=-K  \frac{\partial \zeta_{\mathbf{h}^*}}{\partial \mathbf{n}} \quad
\mbox{on}\, C^k \times ]0,T[,\ \mbox{for}\, k=1,\ldots,N_{CT},
 \vspace{0.1cm}\\ 
\displaystyle
K \frac{\partial \xi}{\partial \mathbf{n}}=b_1^N\big(
\theta_N-\zeta_{\mathbf{h}^*}-\frac{K}{b_1^N} \frac{\partial \zeta_{\mathbf{h}^*}}{\partial \mathbf{n}}-\xi
\big) \quad \mbox{on}\, \Gamma_N \times ]0,T[, 
\vspace{0.1cm}\\ 
\displaystyle
K \frac{\partial \xi}{\partial \mathbf{n}}=b_1^S\big(
\theta_S-\zeta_{\mathbf{h}^*}-\frac{K}{b_1^S} \frac{\partial \zeta_{\mathbf{h}^*}}{\partial \mathbf{n}}-\xi
\big)  \\
\displaystyle
\qquad +b_2^S \big(T_r^4-|\xi+\zeta_{\mathbf{h}^*}|^3(\xi+\zeta_{\mathbf{h}^*})\big) 
\quad \mbox{on}\, \Gamma_S \times ]0,T[, 
\\ 
\displaystyle
\xi(0)=\theta(0)-\zeta_{\mathbf{h}^*}(0) \quad \mbox{in} \, \Omega,
\end{array}\right.
\end{equation}
with $\zeta_{\mathbf{h}^*} \in W^{1,2,2}(0,T;H^{2}(\Omega),
H^{2}(\Omega))$ defined as in Lemma \ref{lemma2}.
\item $\mathbf{h} \in [L^2(0,T)]^{N_{CT}}$ is such that:
\begin{equation}
h^k(s)=\frac{1}{\mu(C^k)} \int_{C^k} \theta(s) \, d \gamma,  \ \mbox{for}\, k=1,2,\ldots,N_{CT}. 
\end{equation}
\end{itemize}

The following technical results are necessary to prove that the operator $M$ defined in (\ref{puntofijo1}) 
is well defined. The first one corresponds to the existence of solution for problem (\ref{eq:system6}), 
and the second one is related to the existence of solution for problem (\ref{eq:system5}). 

\begin{thm} \label{theo2} Within the framework stablished in Hypothesis 
\ref{hypo1}, given elements $\boldsymbol{\zeta}_{\mathbf{g}}\in W^{1,2,2}(0,T;[H_{\sigma}^{2}(\Omega)]^3,[H_{\sigma}^{2}(\Omega)]^3)$ and 
$\theta^* \in L^2(0,T;L^2(\Omega))$, there exists an element $\mathbf{v} \in \mathbf{W}_2$ such that 
$\mathbf{z}=\mathbf{v}-\boldsymbol{\zeta}_{\mathbf{g}}
\in \widetilde{\mathbf{W}}_2$ is the unique solution of problem (\ref{eq:system6}) in the 
following sense:
\begin{equation} \label{eq:system7}
\begin{array}{r}
\displaystyle
\int_{\Omega} \frac{\partial \mathbf{z}}{\partial t} \cdot \boldsymbol{\eta} \,d\mathbf{x} +
\int_{\Omega} \nabla (\boldsymbol{\zeta}_{\mathbf{g}}+
\mathbf{z}) \mathbf{z} \cdot \boldsymbol{\eta} \, d \mathbf{x} 
 \\
 \displaystyle
+
\int_{\Omega} \nabla \mathbf{z} \boldsymbol{\zeta}_{\mathbf{g}} 
\cdot \boldsymbol{\eta} \, d \mathbf{x} 
+ 2 \nu \int_{\Omega} \epsilon(\mathbf{z}) : \epsilon(\boldsymbol{\eta}) \, 
d \mathbf{x} \\
 \displaystyle +
2 \nu_{tur} \int_{\Omega} \left[\epsilon(\boldsymbol{\zeta}_{\mathbf{g}}+\mathbf{z}):
\epsilon(\boldsymbol{\zeta}_{\mathbf{g}}+\mathbf{z}) \right]^{1/2} \epsilon(\boldsymbol{\zeta}_{\mathbf{g}}+\mathbf{z}):
\epsilon(\boldsymbol{\eta}) \, d \mathbf{x} \\
\displaystyle
=\int_{\Omega} \mathbf{H}_{\mathbf{g}} \cdot \boldsymbol{\eta} \, d \mathbf{x},
\quad \mbox{a.e.} \ t \in ]0,T[, \quad \forall
 \boldsymbol{\eta} \in \widetilde{\mathbf{X}}_2,
\end{array}
\end{equation}
with $\mathbf{z}(0)=\mathbf{v}^0$, 
a.e. $\mathbf{x} \in \Omega$, where:
\begin{equation}
\mathbf{H}_{\mathbf{g}}=\alpha_0 (\theta^*-\theta^0) \mathbf{a}_g-
\frac{\partial \boldsymbol{\zeta}_{\mathbf{g}}}{\partial t}
-\nabla \boldsymbol{\zeta}_{\mathbf{g}} \boldsymbol{\zeta}_{\mathbf{g}} 
\in L^2(0,T;[L^2(\Omega)]^3).
\end{equation}
Besides, we have the following estimates:
\begin{equation} \label{eq:acotapf1}
\begin{array}{r}
\displaystyle
\|\mathbf{z}\|_{L^{\infty}(0,T;[L^2(\Omega)]^3)}+
\|\mathbf{z}\|_{L^2(0,T;[W^{1,2}(\Omega)]^3)}+
\|\mathbf{z}\|_{L^3(0,T;[W^{1,3}(\Omega)]^3)}
\\ \leq  \displaystyle 
C_2(\mathbf{v}^0, \theta^0, \mathbf{g})
 \Big[1+\|\theta^*\|_{L^2(0,T;L^2(\Omega))}\Big]
\end{array}
\end{equation}
\begin{equation} \label{eq:acotapf2}
\begin{array}{r}
\displaystyle
\left\|\frac{\partial \mathbf{z}}{\partial t} \right\|_{L^2(0,T;[L^2(\Omega)]^3)} 
+\|\epsilon(\mathbf{z})\|_{
L^{\infty}(0,T;[L^3(\Omega)]^{3 \times 3})}
+\|\epsilon(\mathbf{z})\|_{
L^{\infty}(0,T;[L^2(\Omega)]^{3 \times 3})} 
 \\ 
\displaystyle \leq 
C_3(\mathbf{v}^0, \theta^0, \mathbf{g})\exp\left(1+\|\theta^*\|^2_{L^2(0,T;L^2(\Omega))}\right)
 \Big[1+\|\theta^*\|_{L^2(0,T;L^2(\Omega))} \Big].
\end{array}
\end{equation}
where $C_2$ and $C_3$ are positive constants that depend continuously 
with respect to the space-time configuration of our computational domain, 
$\mathbf{v}^0$, $\theta^0$ and $\mathbf{g}$.
\end{thm}

\begin{pf} It is a direct consequence of Theorem 8 in \cite{fran7}. In this theorem the authors prove the wellposedness
(existence, uniqueness and regularity of solution) for a modified Navier-Stokes system with non-homogeneous Dirichlet boundary conditions, 
by building a continuous extension of the Dirichlet condition and a Galerking approximation of the corresponding problem 
with homogeneous boundary conditions. \hfill $\blacksquare$
\end{pf}

\begin{thm} \label{theo1} Within the framework stablished in Hypothesis 
\ref{hypo1}, given elements $\zeta_{\mathbf{h}^*} \in W^{1,2,2}(0,T;H^2(\Omega),
H^2(\Omega))$ and $\mathbf{v} \in L^{10/3}(0,T;[L^3_{\sigma}(\Omega)]^3)$, 
there exits an element $\theta \in W_1$ such that $\xi=\theta-
\zeta_{\mathbf{h}^*}\in \widetilde{W}_1$ is the unique solution of 
problem (\ref{eq:system5}) in the following sense:
\begin{equation}\label{eq:system8}
 \begin{array}{r}
\displaystyle \int_{\Omega} \frac{\partial \xi}{\partial t} \eta \, d \mathbf{x} +
\int_{\Omega} \mathbf{v} \cdot \nabla \xi \eta \, d \mathbf{x} + 
K \int_{\Omega} \nabla \xi \cdot \nabla \eta \, d \mathbf{x} +  
b_1^N \int_{\Gamma_N} \xi \eta \, d \gamma 
\\ 
\displaystyle
+b_1^S \int_{\Gamma_S} \xi \eta \, d \gamma
+b_2^S \int_{\Gamma_S} |\xi+\zeta_{\mathbf{h}^*}|^3(\xi+\zeta_{\mathbf{h}^*}) \eta 
\, d \gamma = \int_{\Omega} H_{\mathbf{h}^*}  \eta \, d \mathbf{x}
\\ 
\displaystyle
+\int_{\Gamma_C} g_{\mathbf{h}^*}^{C} \eta \, d \gamma
+b_1^N \int_{\Gamma_N}   g_{\mathbf{h}^*}^{N} \eta \, d \gamma
+b_1^S \int_{\Gamma_S}   g_{\mathbf{h}^*}^{S} \eta \, d \gamma
\\ 
\displaystyle
+b_2^S \int_{\Gamma_S} T_r^4 \eta \, d \gamma
,\quad \mbox{a.e.} \; t \in ]0,T[, \quad 
\forall \eta \in \widetilde{X}_1,
\end{array}
\end{equation}
with $\xi(0)=\theta^0-\zeta_{\mathbf{h}^*}(0)$, a.e. 
$\mathbf{x}\in \Omega$. 

Besides, we have the following estimates:
\begin{equation} \label{eq:acotapf3}
\hspace{-.3cm}  \begin{array}{r}
\displaystyle 
\| \xi \|_{L^{\infty}(0,T;L^2(\Omega))}+
\|\xi\|_{L^2(0,T;\widetilde{X}_1))} +
\|\xi\|_{L^5(0,T;L^5(\Gamma_S))} \\
\displaystyle \leq C_4(\theta^0,\theta_N,\theta_S,T_r)\Big[
1+  \|\mathbf{h}^*\|_{[L^2(0,T)]^{N_{CT}}} \Big],
\end{array}
\end{equation}
\begin{equation} \label{eq:acotapf4}
\begin{array}{r}
\displaystyle 
\left\|
\frac{d \xi}{dt}
\right\|_{L^{5/4}(0,T;\widetilde{X}_1')}  \leq 
C_5(\theta^0,\theta_N,\theta_S,T_r)
\Big[
1+  \\ 
\displaystyle \|\mathbf{v}\|_{L^{10/3}(0,T;[L^3(\Omega)]^3)}^2+
\|\mathbf{h}^*\|^2_{[L^2(0,T)]^{N_{CT}}} \Big],
\end{array}
\end{equation}
where $C_4$ and $C_5$ are positive constants that depend continuously 
with respect to the space-time configuration of our computational domain, 
$\theta^0$, $\theta_N$, $\theta_S$ and $T_r$.
\end{thm}

\begin{pf} It is a straightforward consequence of Theorem \ref{Atheo1} 
that we will prove in Appendix \ref{app}. In our case, we have to choose there $\Gamma_R=\Gamma_S$, 
$\Gamma_L=\Gamma_T$ and $\Gamma_A=\Gamma_C \cup \Gamma_N$.
\hfill $\blacksquare$
\end{pf}

\begin{lem} The operator $M$ defined in (\ref{puntofijo1}) is well defined and compact.
\end{lem}

\begin{pf} Thanks to above Theorems \ref{theo2} and \ref{theo1}, it is straightforward that the operator $M$ 
is well defined. Let us check now its compactness. 

So, given a bounded sequence 
$\{(\theta_n^*,\mathbf{h}_n^*)\}_{n \in \mathbb{N}} \subset L^2(0,T;L^2(\Omega))
\times [L^2(0,T)]^{N_{CT}}$, we have, using estimates (\ref{eq:acotapf1}) and (\ref{eq:acotapf2}), that 
the corresponding sequence $\{\mathbf{z}_n\}_{n \in \mathbb{N}} \subset \widetilde{\mathbf{W}}_2$ 
of solutions for the problem (\ref{eq:system7}) is bounded in $\widetilde{\mathbf{W}}_2$. Then, we have, 
taking subsequences if necessary, that:
\begin{itemize}
\item $\theta_n^*\rightharpoonup \theta$ in $L^2(0,T;L^2(\Omega))$,
\item $\mathbf{h}_n^* \rightharpoonup \mathbf{h}$ in $[L^2(0,T)]^{N_{CT}}$,
\item $\mathbf{z}_n \rightarrow \mathbf{z} $ strongly in $L^p(0,T;[L^q(\Omega)]^3)$, for all $1<p<\infty$, and $2\leq q < \infty$,
\item $\displaystyle \mathbf{z}_n \rightharpoonup \mathbf{z}$ weakly in $L^3(0,T;\widetilde{\mathbf{X}})$,
\item $\displaystyle \frac{d \mathbf{z}_n}{dt} \rightharpoonup \frac{d \mathbf{z}}{dt}$ weakly in $L^2(0,T;[L^2(\Omega)]^3)$,
\item $\displaystyle \nabla \mathbf{z}_n \rightharpoonup^* \nabla \mathbf{z}$ weakly-$*$ in $L^{\infty}(0,T;[L^3(\Omega)]^3)$,
\item $\displaystyle \beta(\epsilon(\boldsymbol{\zeta}_{\mathbf{g}}+\mathbf{z}_n)) \, \epsilon(\boldsymbol{\zeta}_{\mathbf{g}}+
\mathbf{z}_n)\rightharpoonup \beta(\epsilon(\boldsymbol{\zeta}_{\mathbf{g}}+\mathbf{z})) \, 
\epsilon(\boldsymbol{\zeta}_{\mathbf{g}}+ \mathbf{z})$ weakly in $L^{3/2}(0,T;\widetilde{\mathbf{X}}')$,
\end{itemize}
where $\mathbf{z} \in \widetilde{\mathbf{W}}_2$ is the solution of (\ref{eq:system7}) associated to 
$\theta$. The last convergence is a consequence of the monotony of operator (see \cite{fran7} for more 
details)
\begin{equation}
\begin{array}{rcl}
A:\widetilde{\mathbf{X}}_2& \rightarrow &\widetilde{\mathbf{X}}_2'  \\ 
\mathbf{z} & \rightarrow & A(\mathbf{z}),
\end{array}
\end{equation}
where, for any $\boldsymbol{\xi}\in \widetilde{\mathbf{X}}_2$,
\begin{equation}
\langle A(\mathbf{z}),\boldsymbol{\xi} \rangle =\int_{\Omega} 
\beta(\epsilon(\boldsymbol{\zeta}_{\mathbf{g}}+\mathbf{z})) 
\epsilon(\boldsymbol{\zeta}_{\mathbf{g}}+\mathbf{z}):\epsilon(\boldsymbol{\xi}) \, 
d \mathbf{x}.
\end{equation}
Now, using the results proved in Lemma \ref{lemma2}, we have that the corresponding 
sequence $\{\zeta_{\mathbf{h}_n^*}\}$ 
converges to $\boldsymbol{\zeta}_{\mathbf{h}}$ strongly in $W^{1,2,2}(0,T;H^2(\Omega),H^2(\Omega))$, 
and, consequently,
\begin{itemize}
\item $H_{\mathbf{h}_n^*} \rightarrow H_{\mathbf{h}}$ in $L^2(0,T;L^2(\Omega))$,
\item $g^C_{\mathbf{h}_n^*} \rightarrow g^C_{\mathbf{h}}$ in $L^2(0,T;L^2(\Gamma_C))$,
\item $g^N_{\mathbf{h}_n^*} \rightarrow g^N_{\mathbf{h}}$ in $L^2(0,T;L^2(\Gamma_N))$, 
\item $g^S_{\mathbf{h}_n^*} \rightarrow g^S_{\mathbf{h}}$ in $L^2(0,T;L^2(\Gamma_S))$.
\end{itemize}
Finally, thanks to estimates (\ref{eq:acotapf3}) and (\ref{eq:acotapf4}), 
the corresponding sequence $\{\xi_n\}_{n \in \mathbb{N}}$ $\subset \widetilde{W}_1$ is 
bounded. Thus, taking subsequences if necessary, we have 
thanks to a straightforward adaptation of Lemma \ref{lemma4} 
in Appendix \ref{app}, that 
\begin{itemize}
\item $\displaystyle \xi_n \rightharpoonup \xi$ in $L^2(0,T;\widetilde{X})$, 
\item $\displaystyle \xi_n \rightharpoonup^* \xi$ in $L^{\infty}(0,T;L^2(\Omega))$, 
\item $\displaystyle \xi_n \rightarrow \xi$ in $L^{10/3-\epsilon}(0,T;L^{10/3-\epsilon}(\Omega))$,
\item $\displaystyle \xi_n \rightarrow \xi$ in $L^2(0,T;L^2(\Gamma_C))$,
\item {\color{red} $\displaystyle \xi_n \rightarrow \xi$ in $L^4(0,T;L^4(\Gamma_R))$.}
\end{itemize}
Using the same techniques that we present in Appendix \ref{app} for the 
demonstration of Theorem \ref{Atheo1}, we can pass to the limit in the 
variational formulation of $\xi_n$ and prove that 
$\xi\in \widetilde{W}_2$ is the solution of (\ref{theo1}) associated to 
$\theta$ and $\mathbf{h}$. Thus, we have that 
\begin{itemize}
\item $\displaystyle \theta_n = \xi_n+\zeta_{\mathbf{h}_n^*} \rightarrow \theta=\xi+
\zeta_{\mathbf{h}}$ in $L^2(0,T;L^2(\Omega))$,
\item $\displaystyle \theta_n = \xi_n+\zeta_{\mathbf{h}_n^*} \rightarrow \theta=\xi+
\zeta_{\mathbf{h}}$ in $L^2(0,T;L^2(\partial \Omega))$,
\end{itemize}
and, consequently, $M(\theta_n^*,\mathbf{h}_n^*)\rightarrow M(\theta,\mathbf{h})$ in $ L^2(0,T;L^2(\Omega)) 
\times [L^2(0,T)]^{N_{CT}}$,
which concludes the proof. \hfill $\blacksquare$
\end{pf}

\begin{thm} Given positive constants $\widehat{C}^1$ and $\widehat{C}^2$, there exist 
coefficients and data small enough such that the operator $M$ defined in (\ref{puntofijo1}) has a fixed point
 in the space $\{(\theta,\mathbf{h}) \in L^2(0,T;L^2(\Omega))\times [L^2(0,T)]^{N_{CT}}:\;
\|\theta\|_{L^2(0,T;L^2(\Omega))}  \leq \widehat{C}_1,\; \|\mathbf{h}\|_{[L^2(0,T)]^{N_{CT}}} \leq 
\widehat{C}_2\}$. Moreover, the corresponding $(\theta,\mathbf{v} )\in W_1 \times \mathbf{W}_2$ is 
a solution for the system (\ref{eq:system1}) and (\ref{eq:system2}) in the sense of 
Definition \ref{defsol}.
\end{thm}

\begin{pf} The existence is a direct consequence of the Schauder fixed point Theorem. Given an 
element $(\theta^*,\mathbf{h}^*) \in L^2(0,T;L^2(\Omega))\times [L^2(0,T)]^{N_{CT}}$, 
we have, thanks to (\ref{eq:acotapf0}), (\ref{eq:acotapf1}), (\ref{eq:acotapf2}),  
(\ref{eq:acotapf3}) and (\ref{eq:acotapf4}), the following estimates for $(\theta,\mathbf{h}) =M(\theta^*,\mathbf{h}^*)$:
\begin{equation}\label{eq:acotapf5}
\begin{array}{l}
\displaystyle 
\|\mathbf{v}\|_{\mathbf{W}_2} \leq   
\displaystyle 
C_6(\mathbf{v}^0,\theta^0,\mathbf{g}) 
\exp(1+\|\theta^*\|_{L^2(0,T;L^2(\Omega))}^2 )
\left[
1+\|\theta^*\|_{L^2(0,T;L^2(\Omega))} 
\right] ,  \\ 
\displaystyle 
\|\theta\|_{W_1} \leq  \displaystyle 
C_7(\theta^0,\theta_N,\theta_S,T_r) \left[
1+\|\mathbf{v}\|_{\mathbf{W}_2}^2 +\|\mathbf{h}^*\|_{[L^2(0,T)]^{N_{CT}}}^2
\right],\\
\displaystyle 
\|\mathbf{h}\|_{[L^2(0,T)]^{N_{CT}}}  \leq 
\displaystyle 
C_8(\theta^0,\theta_N,\theta_S,T_r) \left[
1+\|\mathbf{v}\|_{\mathbf{W}_2}^2 +\|\mathbf{h}^*\|_{[L^2(0,T)]^{N_{CT}}}^2
\right],
\end{array}
\end{equation}
where $C_6$, $C_7$, $C_8$ are positive constants that depend continuously on the coefficients and data. 
If we take the first inequality to the second and third ones, we obtain that:
\begin{displaymath}
\begin{array}{rcl}
\|\theta\|_{W_1} &\leq & \displaystyle 
C_7(\theta^0,\theta_N,\theta_S,T_r) \bigg[
1+\|\mathbf{h}^*\|_{[L^2(0,T)]^{N_{CT}}}^2
\\ && \displaystyle +
C_{9}(\mathbf{v}^0,\theta^0,\mathbf{g})
\exp(1+\|\theta^*\|_{L^2(0,T;L^2(\Omega))}^2 )
\big[
1+\|\theta^*\|_{L^2(0,T;L^2(\Omega))}^2
\big]
\bigg],  \\
\displaystyle 
\|\mathbf{h}\|_{[L^2(0,T)]^{N_{CT}}} & \leq & 
\displaystyle 
C_8(\theta^0,\theta_N,\theta_S,T_r) \bigg[
1+\|\mathbf{h}^*\|_{[L^2(0,T)]^{N_{CT}}}^2 
\\ && \displaystyle +
C_{9}(\mathbf{v}^0,\theta^0,\mathbf{g})
\exp(1+\|\theta^*\|_{L^2(0,T;L^2(\Omega))}^2 )
\big[
1+\|\theta^*\|_{L^2(0,T;L^2(\Omega))}^2
\big]
\bigg].
\end{array}
\end{displaymath}
So, if we suppose that $\|\theta^*\|_{L^2(0,T;L^2(\Omega))} \leq \widehat{C}_1$ and 
$\|\mathbf{h}^*\|_{[L^2(0,T)]^{N_{CT}}}^2 \leq \widehat{C}_2$, we have that
\begin{displaymath}
\begin{array}{l}
\|\theta\|_{W_1} \leq  \displaystyle 
C_7(\theta^0,\theta_N,\theta_S,T_r) \bigg[
1+\widehat{C}_2^2
+
C_{9}(\mathbf{v}^0,\theta^0,\mathbf{g})
\exp(1+\widehat{C}_1^2 )
\big[
1+\widehat{C}_1^2
\big]
\bigg],  \\
\displaystyle 
\|\mathbf{h}\|_{[L^2(0,T)]^{N_{CT}}}  \leq  
\displaystyle 
C_8(\theta^0,\theta_N,\theta_S,T_r) \bigg[
1+\widehat{C}_2^2 +
C_{9}(\mathbf{v}^0,\theta^0,\mathbf{g})
\exp(1+\widehat{C}_1^2 )
\big[
1+\widehat{C}_1^2
\big]
\bigg].
\end{array}
\end{displaymath}
Thus, we are led to solve the following inequality:
\begin{equation}
\begin{array}{l}
 \displaystyle 
C_7(\theta^0,\theta_N,\theta_S,T_r) \bigg[
1+\widehat{C}_2^2
+
C_{9}(\mathbf{v}^0,\theta^0,\mathbf{g})
\exp(1+\widehat{C}_1^2 )
\big[
1+\widehat{C}_1^2
\big]
\bigg] \leq \widehat{C}_1, \vspace{0.1cm} \\ 
\displaystyle 
C_8(\theta^0,\theta_N,\theta_S,T_r) \bigg[
1+\widehat{C}_2^2 +
C_{9}(\mathbf{v}^0,\theta^0,\mathbf{g})
\exp(1+\widehat{C}_1^2 )
\big[
1+\widehat{C}_1^2
\big]
\bigg] \leq \widehat{C}_2.
\end{array}
\end{equation}
However, it is obvious that, given $\widehat{C}_1$ and $\widehat{C}_2$, we can consider 
small enough data $\mathbf{v}^0$, $\mathbf{g}$,  $\theta^0$, 
$\theta_N$, $\theta_S$ and $T_r$,
such that 
\begin{eqnarray}  \label{eq:sch1}
&&C_{9}(\mathbf{v}^0,\theta^0,\mathbf{g}) \leq \frac{1}{\exp(1+\widehat{C}_1^2 )\big[1+\widehat{C}_1^2\big]},  \\
\displaystyle 
&&C_7(\theta^0,\theta_N,\theta_S,T_r)\leq  \displaystyle
\frac{\widehat{C}_1}{2+\widehat{C}_2}, \label{eq:sch4}  \\
\displaystyle
&&C_8(\theta^0,\theta_N,\theta_S,T_r) \leq  \displaystyle
\frac{\widehat{C}_2}{2+\widehat{C}_2}. \label{eq:sch2}
\end{eqnarray}
Then, choosing suitable coefficients and data that verify (\ref{eq:sch1})-(\ref{eq:sch2}), 
we have that $M$ maps elements of the set $\{(\theta^*,\mathbf{h}^*) \in
L^2(0,T;L^2(\Omega))
\times [L^2(0,T)]^{N_{CT}} : \; \|\theta^*\|_{L^2(0,T;L^2(\Omega))}\leq 
\widehat{C}_1,\; \|\mathbf{h}^*\|_{[L^2(0,T)]^{N_{CT}}} 
\leq \widehat{C}_2\}$ into itself. Thus, thanks to Schauder fixed point Theorem, 
there exists a fixed point $(\theta,\mathbf{h})$ of operator $M$, such that the corresponding 
$(\theta,\mathbf{v} )$ is a solution of the coupled system (\ref{eq:system1}) 
and (\ref{eq:system2}) in the sense of Definition \ref{defsol}.  \hfill $\blacksquare$
\end{pf}

\section{Numerical resolution}

Once proved in above section that the coupled system (\ref{eq:system1}) and (\ref{eq:system2}) admits
a solution, we will introduce here a full numerical algorithm in order to compute it,
and show several computational test for a realistic example.

We must recall here that our main aim is related to understanding which is the best strategy for 
reducing the water temperature in the upper layers. In order to achieve this objective, 
and for the sake of completeness, we will consider an algorithm able to deal with more general states 
than those we have presented in previous mathematical analysis of the problem. 
In particular, in the numerical resolution proposed here we will also take into account 
the possibility that $g^k(t)$ takes non-positive values. To be exact, if $g^k(t)>0$ we will say that 
the pump $k$ is turbinating (water enters by the collector $C^k$ and is turbinated by the corresponding pipeline to 
the injector $T^k$), and if $g^k(t)<0$ we will say that the pump $k$ is pumping (water 
enters by injector $T^k$ and is pumped to the collector $C^k$). As it is evident, the situation $g^k(t)=0$ 
corresponds to the case in which the pump $k$ is off. In addition, we will also suppose that the 
parameter $\epsilon$ used in the definition (\ref{eq:defphi}) tends to cero, that is, the mean 
temperature in the injectors is equal to the mean temperature in the collectors.
It is essential emphasizing here that it is also possible to perform a similar mathematical analysis for the general case
(with the obvious embarrassing notations), under the only assumption of the existence of a partition of the time interval 
verifying that the groups do not change their state within any element of the partition (as can be seen in following section
devoted to the numerical examples).

\subsection{Space-time discretization}

For the discretization of the problem, let us consider a regular partition $0=t_0<t_1<\ldots<t_N=T$ of the time interval $[0,T]$ such that
$t_{n+1}-t_n=\Delta t=\frac{1}{\alpha}$, $\forall n =0,\ldots,N-1$, and a family of meshes $\tau_h$
for the domain $\Omega$ with characteristic size $h$.
Associated to this family of meshes, we also consider three compatible finite element spaces $Z_h$, $\mathbf{W}_h$ and $M_h$ 
corresponding, respectively, to the water temperature, velocity and the pressure of water. 
From the computational viewpoint, for the generation of the 
mesh associated to the domain and for the numerical resolution of the system, we propose the use of FreeFem++ \cite{HECHT1}. 
Finally, we have employed an Uzawa algorithm \cite{girault1} 
for computing the solution of the Stokes problems that appears after the 
discretization, and a fixed point algorithm for solving the nonlinearities.

So, we consider the following space-time discretization for system 
(\ref{eq:system1}) and (\ref{eq:system2}):
\begin{enumerate}

\item \textit{Dirichlet condition for the water velocity}: We consider the following approximation 
for a function $\mathbf{g}=(g^1,\ldots,g^{N_{CT}}) \in [H^1(0,T)]^{N_{CT}}$, 
\begin{equation}
g^k(x)=\sum_{n=0}^N g^{k,n} e^n(x),
\end{equation}
where, for all $n=0,\ldots,N$, $e^n \in \mathcal{C}([0,T])$ is such that $e_{|_{[t^{n-1},t_n]}} \in 
\mathcal{P}_1([t^{n-1},t^n])$, $n=1,\ldots,N$ and $e^n(t^k)=\delta_{k,n}$, $k,n=0,\ldots,N$. 
It is well known that the linear closure of the functions of the basis $\{e^n\}_{n=0}^N$ is a vector 
subspace of $H^1(0,T)$, so if we suppose that $g^{k,0}=0$, $k=1,\ldots,N_{CT}$, we can 
consider the following coordinate vector in the basis of the corresponding subspace of 
$[H^1(0,T)]^{N_{CT}}$:
\begin{equation} 
\mathbf{g}=(\underbrace{g^{1,1},g^{2,1},\ldots,g^{N_{CT},1}}_{\mathbf{g}^1},\ldots,
\underbrace{g^{1,N},g^{2,N},\ldots,g^{N_{CT},N}}_{\mathbf{g}^{N}}) \in 
\mathbb{R}^{N\times N_{CT}},
\end{equation}
with $g^{k,n} \in [-M,M]$, $k=1,\ldots,N_{CT}$, $n=1,\ldots,N$, with $M>0$ a 
technical bound related to mechanical characteristics of pumps.

\item \textit{Coupling of temperature in collectors and injectors}:  We denote by $\theta^n\in Z_h$ 
the water temperature at time step $n=0,\ldots,N$. Then, we can consider the following approximation 
in the case of $g^{k,n}>0$, $k=1,\ldots,N_{CT}$, $n=1,\ldots,N$, for (\ref{eq:gammatheta}) by 
functions $\gamma_{\theta}^k$, $k=1,\ldots,N_{CT}$, :
\begin{eqnarray*} 
\gamma_{\theta}^k(t)=\frac{1}{\mu(C^k)} \bigg[ \chi_{(-\infty,t_0)} \int_{C^k} \theta^0 d \gamma 
\vspace{0.2cm} \\
+
\sum_{n=1}^{N} \chi_{[t_{n-1},t_n)} \int_{C^k} \theta^{n-1} d \gamma + \chi_{[t_N,\infty)} \int_{C^k} \theta^N d \gamma \bigg]
\end{eqnarray*}
Moreover, if we assume the value $\epsilon=\frac{\Delta t}{2}$ in the definition (\ref{eq:rhoeps}) of function $\rho_{\epsilon}$ we have 
that the support of $\rho_{\Delta t/2}(t^n-\frac{\Delta t}{2}-s)$ is contained in $(t^n-\Delta t, t^n)=(t^{n-1},t^n)$, for all $n=1,\ldots, N$, and then:
\begin{equation} \nonumber
\begin{array}{rcl}
\displaystyle 
\phi_{\theta}^n(\mathbf{x})&=& \displaystyle \sum_{k=1}^{N_{CT}} \varphi^k(\mathbf{x}) 
\int_{-T}^T \rho_{\Delta t/2}(t^n-\frac{\Delta t}{2} -s) \, \gamma_{\theta}^k(s) \, ds  \\
&=&\displaystyle \sum_{k=1}^{N_{CT}} \varphi^k(\mathbf{x}) \int_{t_{n-1}}^{t_n}  \rho_{\Delta t/2}(t^n-\frac{\Delta t}{2} -s) \left[
\frac{1}{\mu(C^k)} \int_{C^k} \theta^{n-1} \, d \gamma\right] ds  \\
&=&\displaystyle \sum_{k=1}^{N_{CT}} \varphi^k(\mathbf{x}) \, \frac{1}{\mu(C^k)} \int_{C^k} \theta^{n-1} \, d \gamma.
\end{array}
\end{equation}
Finally, we approximate each element $\varphi^k$ by the indicator function of the injector $T^k$, $k=1,\ldots, 
N_{CT}$, and each element $\tilde{\varphi}^k$ by the indicator function of the collector $C^k$, $k=1,\ldots,N_{CT}$. 
Thus, the temperature in each injector at time step $t_n$ is the mean temperature in the corresponding collector at time step $t_{n-1}$. The previous approximation is still valid in the general case with obvious modifications. For 
instance, if we suppose that  $g^{k,n}<0$, $k=1,\ldots,N_{CT}$, $n=1,\ldots,N$, then, we consider 
the Dirichlet condition 
in the collectors:
\begin{equation}  \nonumber
\phi_{\theta}^n(\mathbf{x}) = 
\sum_{k=1}^{N_{CT}} \widetilde{\varphi}^k(\mathbf{x}) \, \frac{1}{\mu(T^k)} \int_{T^k} \theta^{n-1} \, d \gamma.
\end{equation}

\item \textit{Water temperature:} Given $\theta^0 \in Z_h$, $\theta^1 \in 
Z_h$ is the solution of:
\begin{equation} 
\begin{array}{r}
\displaystyle
\alpha \int_{\Omega} \theta^1 \eta \, d \mathbf{x} + 
K \int_{\Omega} \nabla \theta^1 \cdot \nabla \eta \, d \mathbf{x} +
b_1^N \int_{\Gamma_N} \theta^1 \eta \, d \gamma+
b_1^S \int_{\Gamma_S} \theta^1 \eta \, d \gamma
 \\ \displaystyle
+
b_2^S \int_{\Gamma_S} |\theta^1|^3 \theta^1 \eta \, d \gamma
=\alpha \int_{\Omega} (\theta^0 \circ X^0) \eta \, d\mathbf{x} +
b_1^N\int_{\Gamma_N} \theta_N^1 \eta \, d \gamma 
 \\ \displaystyle
+
b_1^S \int_{\Gamma_S} \theta_S^1 \eta\, d \gamma
+b_2^S \int_{\Gamma_S} (T_r^1)^4 \eta \, d \gamma, 
\quad \forall \eta \in Z_h,
\end{array}
\end{equation}
where the discrete characteristic $X^0(\mathbf{x})=
\mathbf{x}-\Delta t \, \mathbf{v}^0(\mathbf{x})$.

Then, for each $n=1,\ldots,N$, $\theta^{n+1} \in Z_h$, with
\begin{eqnarray} \label{eq:bc1}
\displaystyle 
\frac{\partial \theta^{n+1}}{\partial \mathbf{n}}_{|_{C_k}} = 0 , \quad 
\theta^{n+1}_{|_{T^k}}=\frac{1}{\mu(C^k)} \int_{C^k} \theta^{n+1} \, d \gamma
&& \textrm{ if $g^{k,n}>0$},  \\
\label{eq:bc2}
\displaystyle 
\theta^{n+1}_{|_{C^k}}= \frac{1}{\mu(T^k)} \int_{T^k} \theta^{n+1} \, d \gamma, \quad 
\frac{\partial \theta^{n+1}}{\partial \mathbf{n}}_{|_{T^k}} = 0 && \textrm{ if $g^{k,n}<0$},\\
\label{eq:bc3}
\displaystyle 
\frac{\partial \theta^{n+1}}{\partial \mathbf{n}}_{|_{C_k}} = 0, \quad  
\frac{\partial \theta^{n+1}}{\partial \mathbf{n}}_{|_{T^k}} = 0 && \textrm{ if $g^{k,n}=0$}, 
\end{eqnarray}
for all $k=1,\ldots, N_{CT}$, is the solution of:
\begin{equation} \hspace*{-.5cm}
\begin{array}{r}
\displaystyle
\alpha \int_{\Omega} \theta^{n+1} \eta \, d \mathbf{x} + 
K \int_{\Omega} \nabla \theta^{n+1} \cdot \nabla \eta \, d \mathbf{x} +
b_1^N \int_{\Gamma_N} \theta^{n+1} \eta \, d \gamma
\\ \displaystyle + b_1^S \int_{\Gamma_S} \theta^{n+1} \eta \, d \gamma
+
b_2^S \int_{\Gamma_S} |\theta^{n+1}|^3 \theta^{n+1} \eta \, d \gamma
=\alpha \int_{\Omega} (\theta^n \circ X^n) \eta \, d\mathbf{x} 
\\ \displaystyle+
b_1^N\int_{\Gamma_N} \theta_N^{n+1} \eta \, d \gamma
+b_1^S \int_{\Gamma_S} \theta_S^{n+1} \eta\, d \gamma
+b_2^S \int_{\Gamma_S} (T_r^{n+1})^4 \eta \, d \gamma, 
\quad \forall \eta \in Z_h^n,
\end{array}
\end{equation}
where the discrete characteristic $X^n(\mathbf{x})=
\mathbf{x}-\Delta t \, \mathbf{v}^n(\mathbf{x})$, for $n=1,\dots, N$,
and where the functional space $Z_h^n$ is given by:
\begin{equation} \label{eq:zh}
\{\eta \in Z_h: \eta_{|_{
\left(\cup_{k=1}^{N_{CT}} \frac{1}{2}(1-\textrm{sign}(g^{k,n}))C^k\right)|
\textrm{sign}(g^{k,n})|
\cup
\left(\cup_{k=1}^{N_{CT}}\frac{1}{2}(1+\textrm{sign}(g^{k,n}))T^k\right)|
\textrm{sign}(g^{k,n})|}}=0 \},
\end{equation}
with $\textrm{sign}(y)$ denoting the sign function:
\begin{equation}
\textrm{sign}(y)=\left\{ \begin{array}{rl}
1  & \mbox{ if } y > 0, \\
-1 & \mbox{ if } y < 0, \\
0 &  \mbox{ if } y = 0.
\end{array} \right. \nonumber
\end{equation}
For instance, in the case of $g^{k,n}>0$ for all $k=1,\ldots,N_{CT}$, $n=1,\ldots,N$, 
then $\left(\cup_{k=1}^{N_{CT}} \frac{1}{2}(1-\textrm{sign}(g^{k,n}))C^k\right)|
\textrm{sign}(g^{k,n})|$ $\cup
\left(\cup_{k=1}^{N_{CT}}\frac{1}{2}(1+\textrm{sign}(g^{k,n}))T^k\right)$ $= 
\cup_{k=1}^{N_{CT}} T^k$, so
\begin{equation}
Z_h^n=\{\eta \in Z_h:\; \eta_{|_{\cup_{k=1}^{N_{CT}} T^k}}=0\}, \quad n=1,\ldots,N.
\end{equation} 
That is, we are considering the Dirichlet condition on the injectors. In the oposite case, 
$g^{k,n}<0$ for all $k=1,\ldots,N_{CT}$, $n=1,\ldots,N$, we have that
\begin{equation}
Z_h^n=\{\eta \in Z_h:\; \eta_{|_{\cup_{k=1}^{N_{CT}} C^k}}=0\},\quad n=1,\ldots,N,
\end{equation} 
and we are considering the Dirichlet condition on the collectors. In the general case, 
we can have alternating Dirichlet conditions on the collectors and injectors, 
so the definition for the space $Z_h^n$, $n=1,\ldots,N$, given in (\ref{eq:zh}) 
covers all the possibilities.

\item \textit{Water velocity and pressure:} Given $\mathbf{v}^0 \in \mathbf{V}_h$, 
for each $n=0,1,\ldots,N-1$, the pair velocity/pressure $(\mathbf{v}^{n+1},
p^{n+1}) \in \mathbf{W}_h \times M_h$, with:
\begin{equation} 
\begin{array}{r}
\displaystyle
\mathbf{v}^{n+1}_{|{T^k}}=-\frac{g^{k,n+1}}{\mu(T^k)} \,\mathbf{n},
\quad
\mathbf{v}^{n+1}_{|{C^k}}=\frac{g^{k,n+1}}{\mu(C^k)} \,\mathbf{n}, \quad \forall k=1,\ldots,N_{CT},
\end{array}
\end{equation}
is the solution of:
\begin{equation} \label{eq:hidrodiscreto}
\begin{array}{r}
\displaystyle
\alpha \int_{\Omega} \mathbf{v}^{n+1} \cdot \mathbf{z}\, d \mathbf{x} +
2 \nu \int_{\Omega} \epsilon(\mathbf{v}^{n+1}):\epsilon(\mathbf{z}) \, d \mathbf{x} \\
\displaystyle 
+ 2 \nu_{tur}
\int_{\Omega}[\epsilon(\mathbf{v}^{n+1}):\epsilon(\mathbf{v}^{n+1})]^{1/2}\epsilon(\mathbf{v}^{n+1}):
\epsilon(\mathbf{z}) \, d \mathbf{x} \\
\displaystyle
-\int_{\Omega} p^{n+1} \nabla \cdot \mathbf{z}\, d \mathbf{x}
-\int_{\Omega} \nabla \cdot \mathbf{v}^{n+1} q \, d \mathbf{x} 
=\alpha \int_{\Omega} (\mathbf{v}^n \circ X^n) \cdot \mathbf{z}\, d \mathbf{x} \\
\displaystyle
+\int_{\Omega} \alpha_0(\theta^{n+1}-\theta^0) \mathbf{a}_g 
\cdot \mathbf{z}\, d \mathbf{x},
\quad
\forall \mathbf{z} \in \mathbf{V}_h, \quad \forall q \in M_h,
\end{array}
\end{equation}
where the functional space $\mathbf{V}_h=\{\mathbf{z}\in \mathbf{W}_h:\; \mathbf{z}_{|{\cup_{k=1}^{N_{CT}} (T^k \cup C^k) }}=\mathbf{0}\}$.
\end{enumerate}

\begin{rem} It is worthwhile noting here that in above scheme we have to compute one additional time step 
for the water temperature. This shift is motivated by the dependence diagram shown in Fig.~\ref{dscheme}.
We observe that, due to the time discretization proposed here, the Dirichlet boundary condition 
for the hydrodynamic model begins to have influence from the second step of time for 
the water temperature. In the first time step $n=1$, water circulation does not affect the 
water temperature, so for this first time step we consider  
$\Gamma_N=\partial \Omega \setminus \Gamma_S$. However, for $n \geq 2$, we 
impose at each of the $N_{CT}$ collector/injector pairs one of the 
boundary conditions (\ref{eq:bc1}), (\ref{eq:bc2}) or (\ref{eq:bc3}), depending on the 
sign of the Dirichlet condition $g^{k,n-1}$ in previous time step. 

\begin{figure}[!ht]
\[
\xymatrix{  
                              & \mathbf{v}^0 \ar[dr]  \ar[d] & \theta^0 \ar[d] \ar[dl] \\ 
\mathbf{g}^1 \ar[r]  & \mathbf{v}^1 \ar[dr] \ar[d] & \theta^1 \ar[d] \ar[dl] \\ 
\mathbf{g}^2 \ar[r]  & \mathbf{v}^2 \ar[dr] \ar[d] & \theta^2 \ar[d] \ar[dl] \\ 
\mathbf{g}^3 \ar[r]  & \mathbf{v}^3 \ar@{.>}[d]  & \theta^3 \ar@{.>}[d]\\
\mathbf{g}^N \ar[r] & \mathbf{v}^N \ar[dr]         & \theta^N \ar[d] \\ 
                               &                                        &  \theta^{N+1} }
\]
\caption{Dependence scheme for the discretized variables.}
\label{dscheme}
\end{figure}
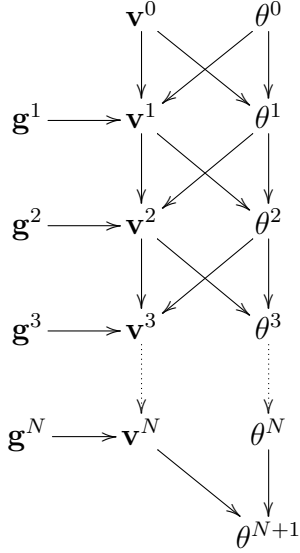

\end{rem}

\subsection{Numerical results}

\begin{figure}[ht]
\centering
\includegraphics[width=.95\linewidth]{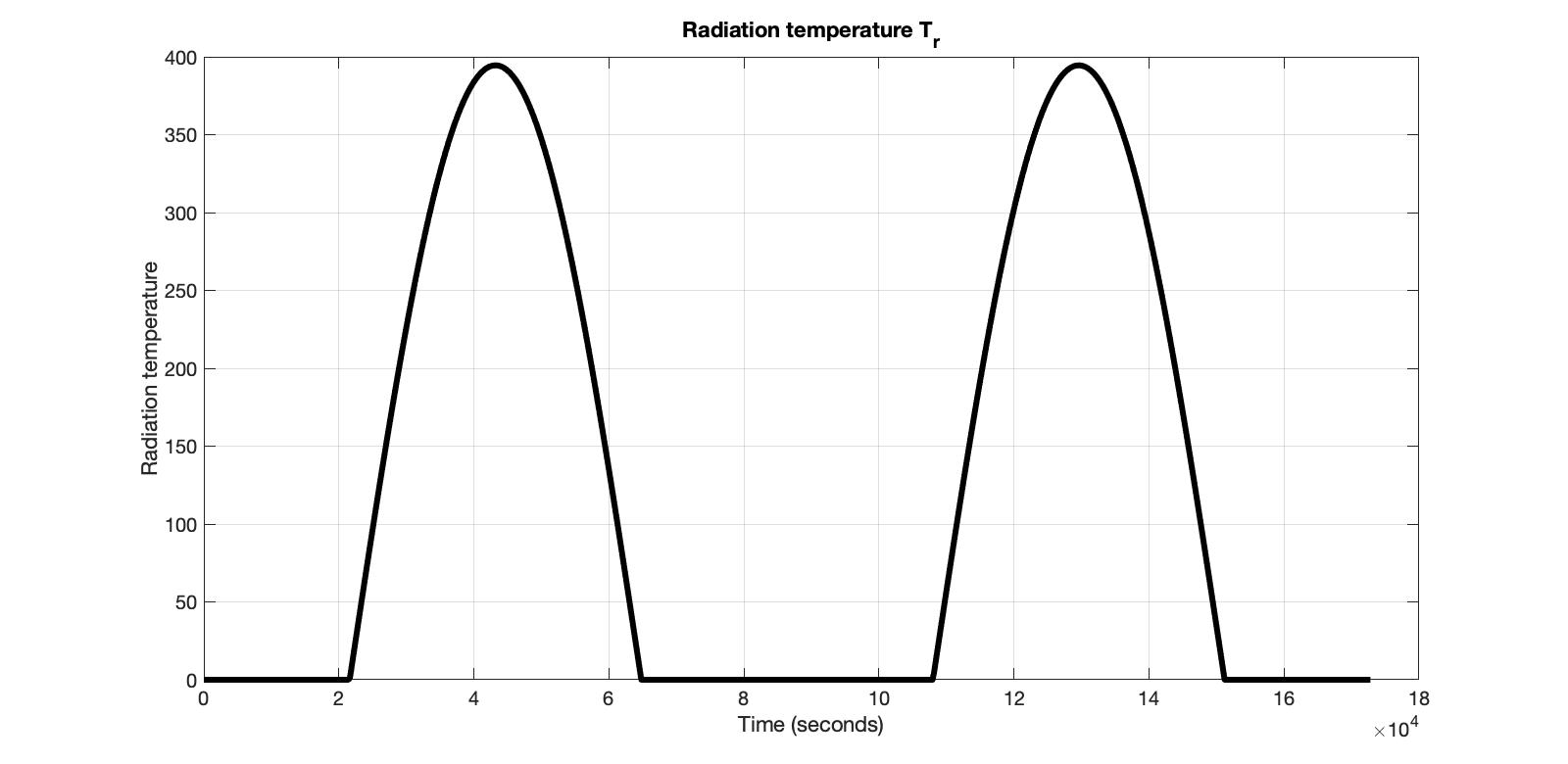}
\caption{Profile of the radiation temperature $T_r$ for the whole time interval of $T=17.28 \ 10^4$ seconds (2 days).}
\label{figure2}
\end{figure}

This final subsection is devoted to present some numerical results that we have 
obtained using realistic coefficients and data. Nevertheless, for the sake of clarity and
comprehensibility, we will show results for a simplified 2D domain. For this 
purpose, we have considered a rectangular domain $\Omega=[0,16] 
\times [0,19]$ (measured in meters), corresponding a reservoir, in which we have distributed 
$N_{CT}=4$ collector/injector pairs with a symmetrical configuration 
similar to that shown in Fig.~\ref{figure1}. For the time discretization 
we have chosen a time step of $\Delta t=1800$ seconds with $N=96$ time steps
(which represents a time period of 2 days), and for the 
space discretization we have used a regular mesh formed by triangles of 
characteristic size $h=0.5$ meters (corresponding to $1599$ vertices). Finally, the 
finite element spaces employed for space discretizations have been 
the Taylor-Hood element $\mathcal{P}_2/\mathcal{P}_1$ for the hydrodynamic 
model, and the Lagrange $\mathcal{P}_2$ element for the water temperature. 
In Fig.~\ref{figure2} we can observe the evolution of the radiation 
temperature:
\begin{equation}
T_r=\left(\frac{1}{\varepsilon\, \sigma_B} \left[ 
(1-a) R_{sw,net}+R_{lw,down}\right] \right)^{1/4}
\end{equation}
along the whole period of 2 days ($17.28 \ 10^4$ seconds), considering 
$a=0.1$, $R_{sw,net}=1000\, {\rm W} \, {\rm m}^{-2}$, $R_{lw,down}=350\, 
{\rm W} \, {\rm m}^{-2}$ (typical values in mediterranean countries during the 
summer) and multiplied by a sinusoidal function in order to simulate 
the effects of day and night. The parameters used for the numerical resolution of the coupled system can be seen in Table \ref{Table}. 

\begin{table}
\centering
\begin{tabular}{|l|l|l|}
\hline Parameters & Values & Units \\ \hline  %\hline
% & & \\
$\nu$ & $1.3 \ 10^{-3}$ & ${\rm m}^2 \, {\rm s}^{-1}$ \\ %\hline
$\nu_{tur}$ & $5.0 \ 10^{-2}$ & ${\rm m}^2$ \\ %\hline
$K$ & $1.4\ 10^{-5}$ & ${\rm m}^2 \, {\rm s}^{-1}$ \\
$h^N$, $h^S$ & $3.0 \ 10^2$ & ${\rm W} \, {\rm m}^{-2} \, {\rm K}^{-1}$ \\
$\rho$ & $9.9 \ 10^2$ & ${\rm g} \, {\rm m}^{-3}$ \\
$c_p$ & $4.2$ & ${\rm W} \, {\rm s} \, {\rm g}^{-1} \, {\rm K}^{-1}$ \\
$\theta^0$ & $283.0$ & K \\ 
$\theta_S$ & $286.0$ & K \\  
$\theta_N$ & $283.0$ & K \\ 
$\alpha_0$ & $8.7\ 10^{-7}$ & ${\rm K}^{-1}$ \\  \hline
\end{tabular}
\caption{Physical parameters for the numerical example.} \label{Table} %\vspace*{-10pt}
\end{table}

In order to analyze the influence of water artificial circulation in the thermal 
behavior of top $1.5$ meters from water upper layer we have solved the problem in 
five different scenarios:
\begin{enumerate}
\item NNNN: In this configuration we take $g^{k,n}=0$, for all $k=1,\ldots,N_{CT}$ 
and $n=1,\ldots,N$ (reference configuration with all the groups off).
\item TTTT: In this configuration we take $g^{k,n}=2.0 \ 10^{-3}$, 
for all $k=1,\ldots,N_{CT}$ and $n=1,\ldots,N$ (all the groups are turbinating). 
\item PPPP: In this configuration we take $g^{k,n}=-2.0\ 10^{-3}$, 
for all $k=1,\ldots,N_{CT}$ and $n=1,\ldots,N$ (all the groups are pumping). 
\item TPTP: In this configuration we take $g^{1,n}=g^{3,n}=2.0\ 10^{-3}$ 
and $g^{2,n}=g^{4,n}=-2.0\ 10^{-3}$, for all $k=1,\ldots,N_{CT}$ and 
$n=1,\ldots,N$ (groups $1$ and $3$ are turbinating, and groups $2$ and $4$ are pumping).
\item PTPT: In this configuration we take $g^{1,n}=g^{3,n}=-2.0\ 10^{-3}$ 
and $g^{2,n}=g^{4,n}=2.0\ 10^{-3}$, for all $k=1,\ldots,N_{CT}$ and the groups $2$ and $4$ are turbinating).
\end{enumerate}

\begin{figure}[ht]
\centering
\includegraphics[width=1.0\linewidth]{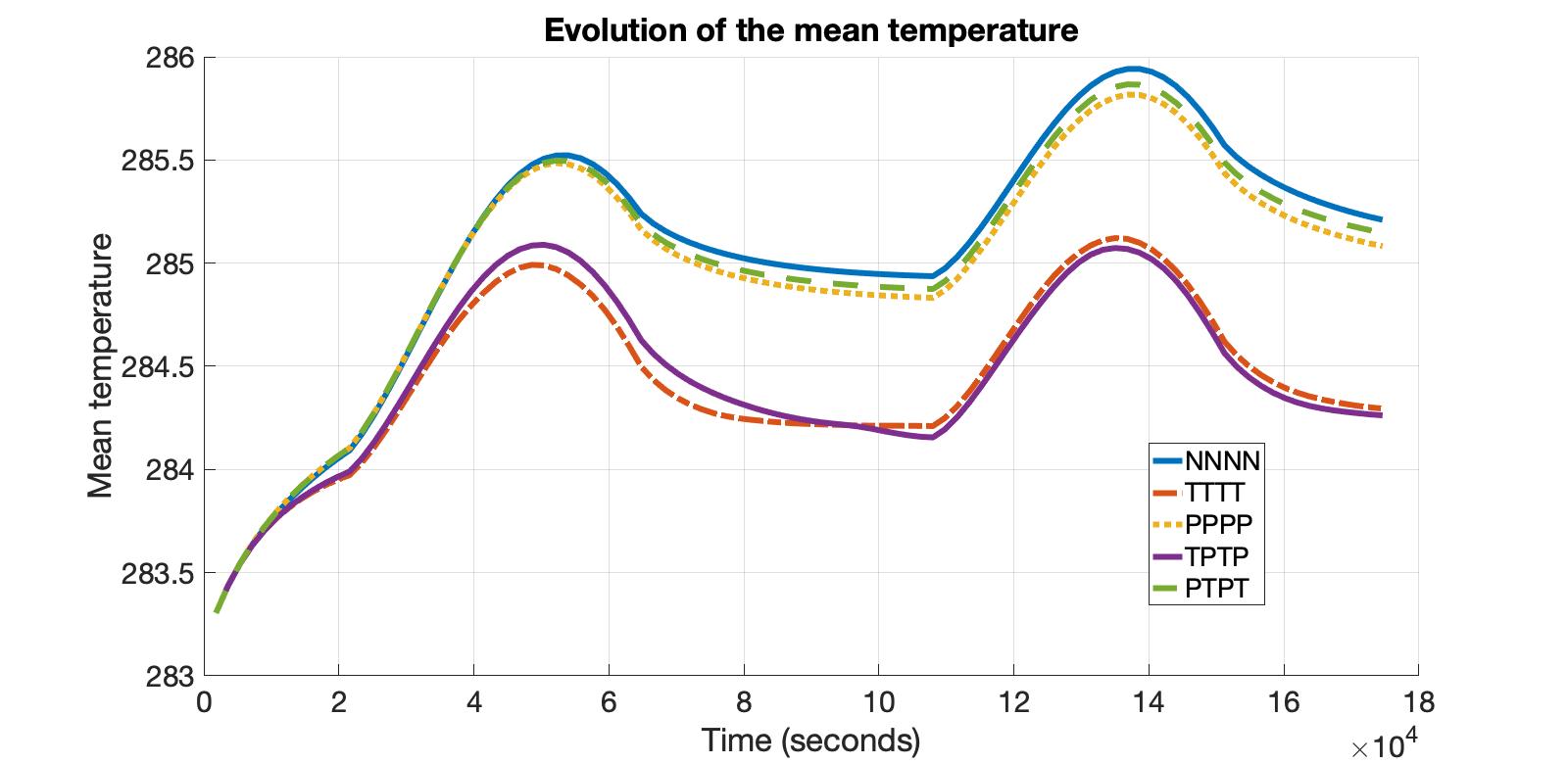}
\caption{Evolution of the mean temperature in the top $1.5$ meters upper layer of the water domain
for the five scenarios under study.}
\label{figure3}
\end{figure}

In Fig.~\ref{figure3} we present the evolution of the mean temperature in 
the top $1.5$ meters upper layer along the whole time interval corresponding to two days. We can clearly distinguish here that the best 
configurations correspond, in a very evident manner, to the second scenario (TTTT) and to the fourth one (TPTP).
Moreover, we can notice how third and fifth scenarios (PPPP and PTPT, respectively) do not improve in a significant way the reference configuration (NNNN).

\begin{figure}[ht]
  \centering
  \includegraphics[width=.44\linewidth]{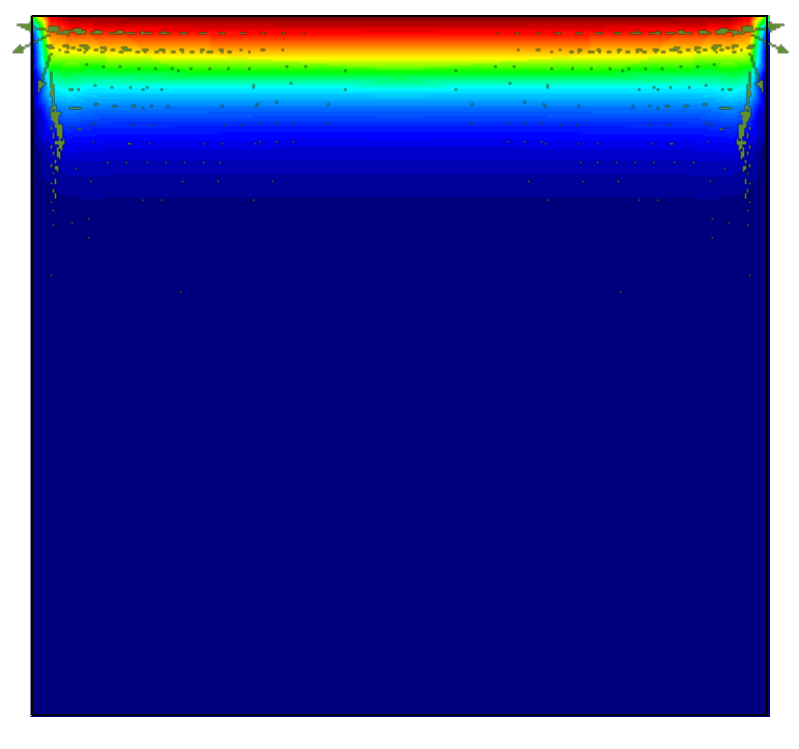}
  \includegraphics[width=.1\linewidth]{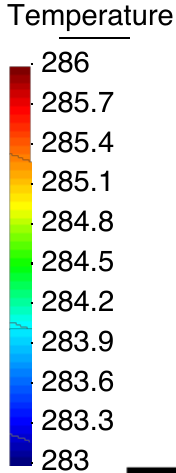}
  \includegraphics[width=.44\linewidth]{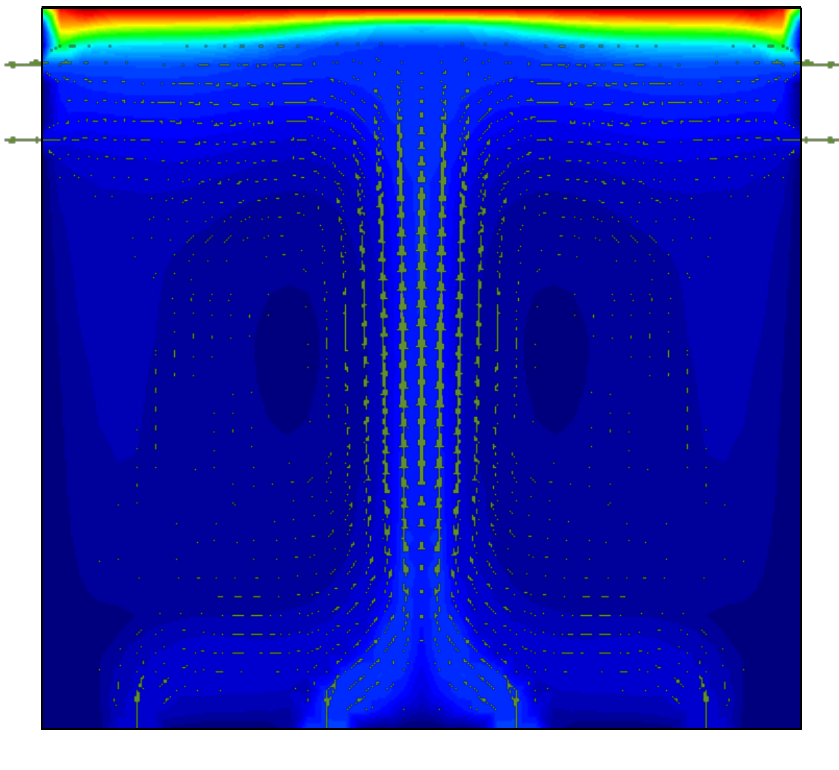}
\caption{Water temperatures and velocities (multiplied by factors 1500/500) in the last time step for 
the worst configuration NNNN (left) and the best configuration TTTT (right).}
\label{figure8}
\end{figure}

\begin{figure}[ht]
  \centering
  \includegraphics[width=.44\linewidth]{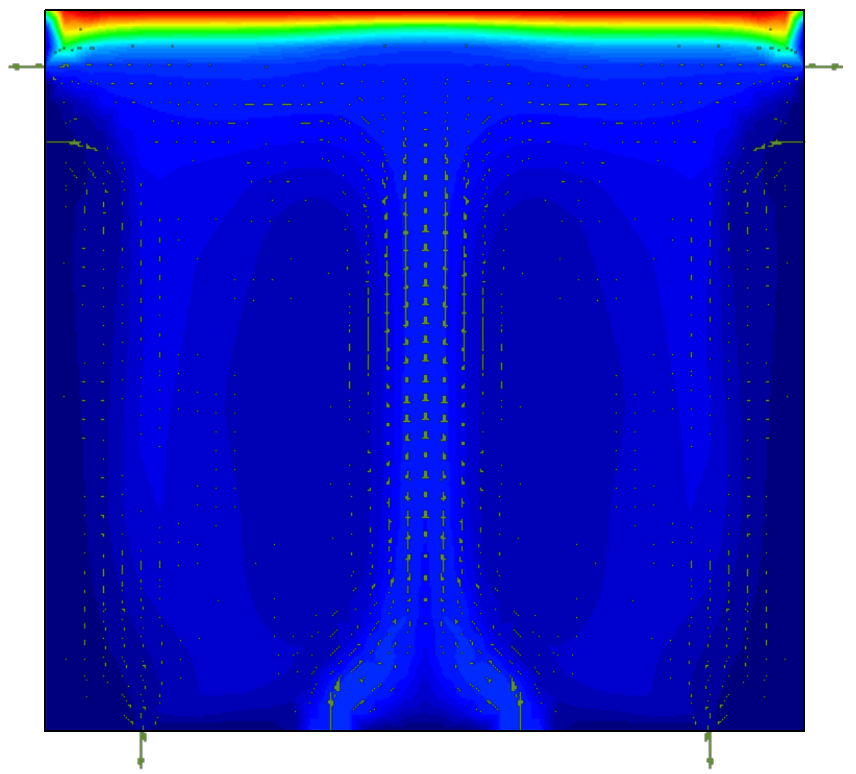}
  \includegraphics[width=.1\linewidth]{image_5.pdf}
  \includegraphics[width=.44\linewidth]{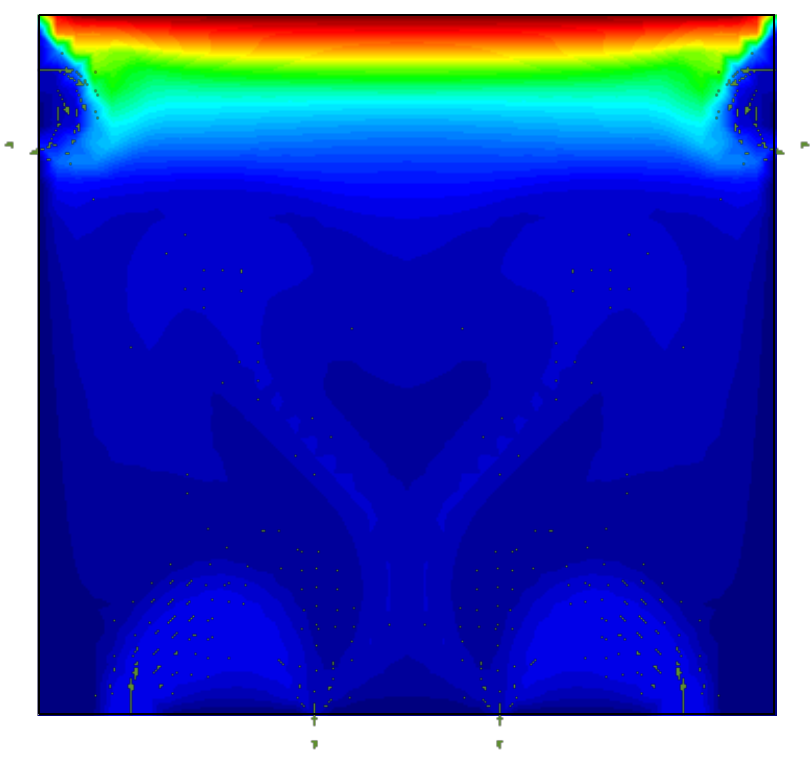}
\caption{Water temperatures and velocities (multiplied by a factor 500) in the last time step for 
the ``good'' configuration TPTP (left) and the ``bad'' configuration PTPT (right).}
\label{figure6}
\end{figure}

Finally, we show in Fig.~\ref{figure8} water temperatures and velocities at last time step for  
NNNN and TTTT configurations and, in Fig.~\ref{figure6}, 
the behavior of water at same time step for configurations TPTP and PTPT. 
(In all of the cases, velocities have been multiplied by an amplifying factor to make their graphic representations more perceptible).
As we can easily notice, the best strategies correspond to evacuating the excess of temperature in the 
upper layers to the bottom layers instead of refrigerating the upper layers with cold water from the bottom ones.

Although we only present here one realistic example of application of our approach to understand
the behaviour of water velocity and temperature in the upper section of the domain, we have developed many other
numerical experiences for different choices of parameters and data (that will not be presented here for the sake of conciseness).
However, from these computational tests we can derive two important consequences: 
For the modified Navier-Stokes equations, the second member -corresponding to the thermic term- shows less influence in the
numerical resuls than the Smagorinsky turbulence term. For the convective heat equation, the radiation term in the nonlinear
boundary condition affects in a significative way the final results.
Finally, for a better resolution of the numerical examples, it would be possible to use, instead of an uniform mesh like the one employed in
previous example, a finer mesh in the neighbourhoods of collectors and injectors. Nevertheless, this approach would mean a significant 
increase in the computational time, already quite high in the current case (especially in the part referring to the resolution of the hydrodynamic problem).

\section*{Acknowledgments}  
The authors thank the funding from project MTM2015-65570-P of MINECO/ FEDER (Spain).

\appendix

\section{A radiation heat transfer problem with nonhomogeneous mixed boundary conditions} \label{app}

In this appendix we mathematically analyze a heat equation with an advective 
term and mixed boundary conditions of diffusive, convective and radiative type. 
The essential difficulty for demonstrating the existence and uniqueness 
of solution lies, in one hand, in the presence of nonlinear boundary 
conditions of radiative type and, in the other hand, the lack of regularity 
of the time derivative of the solutions. This lack of regularity does not 
allow us to take the solution itself as a test function in the variational 
formulation, and we are led to use more refined techniques (similar to 
those employed, for instance, in \cite{fran1}). 

So, we suppose that we have a convex domain $\Omega\subset 
\mathbb{R}^3$, whose boundary can be split into three disjoint, smooth enough parts: $\Gamma_L$, $\Gamma_A$ and 
$\Gamma_R$, with $\partial \Omega=\Gamma_L \cup \Gamma_A \cup \Gamma_R$. We denote 
by $\theta$ the solution of the following initial-boundary value problem:
\begin{equation}\label{eq:Asystem12}
\left\{\begin{array}{l}
\displaystyle \frac{\partial \theta}{\partial t} +\mathbf{v} \cdot \nabla \theta 
-\nabla \cdot (K \nabla \theta) = f \quad \mbox{in}\; \Omega \times ]0,T[, 
 \\  
\theta=\theta_L \quad  \mbox{on}\; \Gamma_L \times ]0,T[, 
\\  
\displaystyle K \nabla \theta \cdot \mathbf{n}=b_1^A (\theta_A-\theta) \quad 
\mbox{on}\; \Gamma_A \times ]0,T[, 
 \\  
\displaystyle K \nabla \theta \cdot \mathbf{n}=b_1^R (\theta_R-\theta)+
b_2^R(\varphi^4-|\theta|^3 \theta) \quad 
\mbox{on}\; \Gamma_R \times ]0,T[, 
\\  
\theta(0)=\theta^0 \quad  \mbox{in}\; \Omega,
\end{array}\right.
\end{equation}
where $K>0 \ ({\rm m}^2 \, {\rm s}^{-1})$ is the thermal diffusivity, 
$b_1^H\geq 0 \ ({\rm m} \, {\rm s}^{-1})$, for $H=A,R$, are the coefficients 
related to convective heat transfer through the boundaries $\Gamma_A$ and 
$\Gamma_R$, and $b_2^R>0\ ({\rm m} \, {\rm s} \, {\rm K}^{-3})$ is the coefficient 
related to radiative heat transfer through the boundary $\Gamma_R$, 
$\theta^0\geq 0 \ ({\rm K})$ is the initial temperature, $\theta_L\geq 0 \ ({\rm K})$ 
is Dirichlet temperature on $\Gamma_L$, $\varphi \geq 0 \ ({\rm K})$ is 
the radiation temperature on $\Gamma_R$, and $\theta_A,\, \theta_R\geq 0 
\ ({\rm K})$ are the temperatures related to convection heat transfer in 
surfaces $\Gamma_A$ and $\Gamma_R$.

\begin{rem} In this work we will suppose that $\Gamma_L$, $\Gamma_A$ and $\Gamma_R$ are nonempty, 
but all the results can be easily extended to the case where $\Gamma_A$ and/or $\Gamma_R$ are empty sets. The only 
drawback is when $\Gamma_L=\emptyset$ because in this case we cannot use Poincare type inequalities, and 
we should apply another type of techniques for obtaining energy estimates in the Galerkin approximation. A related problem 
with $\Gamma_L=\emptyset$ was studied, for instance, in \cite{Metzger1}.
\end{rem}

We consider the following spaces
\begin{equation}
\begin{array}{rcl}
X&=&\{\theta \in H^1(\Omega):\; \theta_{|_{\Gamma_R}}\in L^5(\Gamma_R)\}, \\  
\widetilde{X} &=& \{\theta \in X:\; \theta_{|_{\Gamma_L}}=0\},
\end{array}
\end{equation}
and 
\begin{equation} 
\begin{array}{rcl}
\displaystyle
W&=&\displaystyle
\{\theta \in W^{1,2,5/4}(0,T;X,X'):\\ && \displaystyle
\theta_{|_{\Gamma_R}} \in 
L^5(0,T;L^5(\Gamma_R))\}\cap L^{\infty}(0,T;L^2(\Omega)),  \\ 
\widetilde{W} &=& \displaystyle
\{\theta \in W^{1,2,5/4}(0,T;\widetilde{X},\widetilde{X}'):\\ && \displaystyle
\theta_{|_{\Gamma_R}} \in 
L^5(0,T;L^5(\Gamma_R))\}\cap L^{\infty}(0,T;L^2(\Omega)). 
\end{array}
\end{equation}

\begin{hypo} \label{hypo1n} We will assume the following hypotheses for the coefficients and data:
\begin{enumerate}
\item $\theta^0 \in L^2(\Omega)$.
\item $f \in L^2(0,T;L^2(\Omega))$.
\item $\mathbf{v} \in L^{10/3}(0,T;[L^3_{\sigma}(\Omega)]^3)$.
\item $\theta_A \in L^2(0,T;L^2(\Gamma_A))$.
\item $\theta_R \in L^2(0,T;L^2(\Gamma_R))$.
\item $\varphi \in L^5(0,T;L^5(\Gamma_R))$.
\item $\theta_L(\mathbf{x},t)=g(t) \theta_D(\mathbf{x})$, a.e. $(\mathbf{x},t) \in \Gamma_L \times ]0,T[$, with 
$g \in H^1(0,T)$ and $\theta_D = {\widehat{\theta_D}}_{|_{\Gamma_L}}$, where $\widehat{\theta_D} 
\in H^{3/2}(\partial \Omega)$.  
\end{enumerate}
\end{hypo}
\begin{defn} \label{defsoln} Within the framework established in Hypothesis \ref{hypo1n}, we say that 
an element $\theta \in W$ is a solution of the system (\ref{eq:Asystem12}) 
if there exists $\xi \in \widetilde{W}$ such that:
\begin{itemize}
\item $\theta=\zeta_{D}+\xi$, with $\zeta_{D}=g\, \beta_0(\widehat{\theta_D}) \in W^{1,2,2}(0,T;H^2(\Omega), 
H^2(\Omega))$, where $\beta_0$ is the right inverse of trace operator 
$\gamma_0:H^{2}(\Omega)
\rightarrow H^{3/2}(\partial \Omega)$.
\item $\xi(0)=\theta^0-\zeta_{D}(0)$, a.e. $\mathbf{x} \in \Omega$.
\item $\xi \in \widetilde{W}$ is the solution of the following variational formulation:
\begin{equation}\label{eq:Asystem2}
\begin{array}{r}
\displaystyle \int_{\Omega} \frac{\partial \xi}{\partial t} \eta \, d \mathbf{x} +
\int_{\Omega} \mathbf{v} \cdot \nabla \xi \eta \, d \mathbf{x} + 
K \int_{\Omega} \nabla \xi \cdot \nabla \eta \, d \mathbf{x} 
\\  
\displaystyle +  
b_1^A \int_{\Gamma_A} \xi \eta \, d \gamma 
+b_1^R \int_{\Gamma_R} \xi \eta \, d \gamma
+b_2^R \int_{\Gamma_R} |\xi+\zeta_{D}|^3(\xi+\zeta_{D}) \eta 
\, d \gamma 
\\  
\displaystyle 
= \int_{\Omega} H_{D}  \eta \, d \mathbf{x}
+b_1^A \int_{\Gamma_A}   g_{D}^{A} \eta \, d \gamma
\\  
\displaystyle
+b_1^R \int_{\Gamma_R}   g_{D}^{R} \eta \, d \gamma
+b_2^R \int_{\Gamma_R} \varphi^4 \eta \, d \gamma,
\;
 \mbox{a.e.} \; t \in ]0,T[, \quad 
\forall \eta \in \widetilde{X},
\end{array}
\end{equation}
where some of previous integrals must be understood as duality pairs, and
\begin{equation} \label{eq:Asystem3}
\begin{array}{rcl}
H_{D}&=&\displaystyle -\mathbf{v} \cdot \nabla \zeta_{D}+
\frac{\partial \zeta_{D}}{\partial t} + 
\nabla \cdot (K \nabla \zeta_{D}) \in 
L^2(0,T;L^2(\Omega)), \\ 
\displaystyle
g^A_{D}&=&\displaystyle
\theta_A-\zeta_{D}-\frac{D}{b_1^A} \nabla \zeta_{D} \cdot \mathbf{n}
\in L^2(0,T;L^2(\Gamma_A)), 
\\  
\displaystyle
g^R_{D}&=&\displaystyle
\theta_R-\zeta_{D}-\frac{D}{b_1^R}\nabla \zeta_{D} \cdot \mathbf{n}
\in L^2(0,T;L^2(\Gamma_R)).
\end{array}
\end{equation}
\end{itemize}
\end{defn}

We have the following result that we will prove in the following subsections:

\begin{thm} \label{Atheo1} Within the framework established in 
Hypothesis \ref{hypo1n}, there exists a unique solution $\xi \in \widetilde{W}$ of equation (\ref{eq:Asystem2}) in the sense of Definition 
\ref{defsoln}. Moreover, there exists a constant $C>0$, such that this solution 
satisfies the following inequalities:
\begin{equation} \label{eq:boundcont1}
\begin{array}{r}
\displaystyle 
\| \xi \|^2_{L^{\infty}(0,T;L^2(\Omega))}+
\|\xi\|^2_{L^2(0,T;H^1(\Omega))} +
\|\xi\|^5_{L^5(0,T;L^5(\Gamma_S))}  \\  
\displaystyle \leq
C\Big[ \|\theta^0-\zeta_{D}(0)\|_{L^2(\Omega)}^2+
\|H_D\|^2_{L^2(0,T;L^2(\Omega))}+ \|g_{D}^{A}\|_{L^2(0,T;L^2(\Gamma_A))}^2
 \\  
\displaystyle
+
\|g_{D}^{R}\|_{L^2(0,T;L^2(\Gamma_R)}^2+
\|\zeta_{D}\|_{L^5(0,T;L^5(\Gamma_R))}^5
+\|\varphi\|_{L^5(0,T;L^5(\Gamma_R))}^5 \Big],
\end{array}
\end{equation}
\begin{equation} \label{eq:boundcont2} \hspace{-.5cm}
\begin{array}{r}
\displaystyle 
\left\|
\frac{d \xi}{dt}
\right\|_{L^{5/4}(0,T;\widetilde{X}')} \leq C \Big[
\|H_D\|^2_{L^2(0,T;L^2(\Omega))}
+\|g_{D}^{A}  \|_{L^2(0,T;L^2(\Gamma_A))}
\\  
\displaystyle
+\| g_{D}^{R}  \|_{L^2(0,T;L^2(\Gamma_R))}
+\|\varphi\|^4_{L^5(0,T;L^5(\Gamma_R))}
+\|\mathbf{v}\|_{L^{10/3}(0,T;[L^3(\Omega)]^3)}
\|\xi\|_{L^2(0,T;\widetilde{X})}
\\  
\displaystyle +
\|\xi\|_{L^2(0,T;\widetilde{X})}
+\|\xi\|^4_{L^5(0,T;L^5(\Gamma_R))}
+\|\zeta_{D}\|^4_{L^5(0,T;L^5(\Gamma_R))}\Big].
\end{array}
\end{equation}
\end{thm}

In order to better understand the proof of the previous result, we will divide it 
in five parts: in the first part we will prove some technical results related to
the space where we look for the solution. In the second part we will 
obtain the Galerkin approximation of problem (\ref{eq:Asystem12}). In the 
third part we will analyze the differential equation obtained from the Galerkin 
discretization. In the fourth part we will derive the convergence of the Galerkin 
approximation in a suitable space. Finally, in fifth part we will prove Theorem 
\ref{Atheo1}.

\subsection{Part 1: Some technical results}

We have the following lemma that we will use in following subsections:

{\color{red}
\begin{lem} \label{lemma1n} Let $Z \subset \tilde{X}$ a Banach space 
such that $L^5(\Gamma_R)\subset Z'$. The following inclusion is compact:
\begin{equation} 
\begin{array}{c}
\displaystyle
\{\theta \in W^{1,2,5/4}(0,T;\widetilde{X},Z'):\\
\displaystyle 
\theta_{|_{\Gamma_R}} \in 
L^5(0,T;L^5(\Gamma_R))\}
\subset \subset
L^4(0,T;L^{4}(\Gamma_R))
\end{array}
\end{equation}
\end{lem}

\begin{pf} 
The proof is a direct consequence of Aubin and Lions Lemma 
(see, for instance, Lemma 7.7 of \cite{Roubicek1}), the compactness of 
$\widetilde{X}$ in $L^2(\Gamma_R)$ and the fact that $L^4(0,T;L^4(\Gamma_R))$ is 
is an interpolant between $L^2(0,T;L^2(\Gamma_R))$ and $L^5(0,T;L^5(\Gamma_R))$. 
Indeed, given a bounded sequence $\{\theta_n\}_{n \in \mathbb{N}}$ 
in $\{\theta \in W^{1,2,5/4}(0,T;\widetilde{X},\widetilde{X}'):
\theta_{|_{\Gamma_R}} \in 
L^5(0,T;L^5(\Gamma_R))\}$ we have that there exists a subsequence, 
that we will still denote in the same way, such that 
${\theta_n}_{|_{\Gamma_R}} \rightarrow {\theta}_{|_{\Gamma_R}}$ 
in $L^2(0,T;L^2(\Gamma_R))$ and 
${\theta_n}_{|_{\Gamma_R}} \rightharpoonup 
{\theta}_{|_{\Gamma_R}}$ in $L^5(0,T;L^2(\Gamma_R))$. 
We realize that ${\theta_n}_{|_{\Gamma_R}} \rightarrow \theta_{|_{\Gamma_R}}$ 
in $L^4(0,T;L^4(\Gamma_R))$ 
because:
\begin{equation}
\begin{array}{c}
\displaystyle
\|{\theta_n}_{|_{\Gamma_R}} -{\theta}_{|_{\Gamma_R}}\|_{L^4(0,T;L^4(\Gamma_R))} \\
\displaystyle
\leq \|{\theta_n}_{|_{\Gamma_R}}-{\theta}_{|_{\Gamma_R}}\|^{5/6}_{L^5(0,T;L^5(\Gamma_R))} 
\|{\theta_n}_{|_{\Gamma_R}}-{\theta}_{|_{\Gamma_R}}\|^{1/6}_{L^2(0,T;L^2(\Gamma_R))} \rightarrow 0,
\end{array}
\end{equation} 
when $n \to \infty$. 
\hfill $\blacksquare$ \end{pf}}

\begin{rem} We have that 
$W^{1,2,2}(0,T;H^2(\Omega),H^2(\Omega)) \subset 
\mathcal{C}([0,T];H^2(\Omega))$ (in fact, it is well known that if $p,q \geq 1$ and $V_1 
\subset V_2$ continuously, then $W^{1,p,q}(0,T;$ $V_1,V_2) 
\subset \mathcal{C}([0,T];V_2)$ continuously), and that $\mathcal{C}([0,T];H^2(\Omega))
\subset W$. Then, the sum $\theta_D+\xi$ makes sense in the space $W$. 
Moreover, $\mathbf{v} \cdot \nabla \zeta_D \in L^{10/3}(0,T;L^{2}(\Omega))$. 

We can also consider as a Dirichlet condition 
the restriction to $\Gamma_L$ of one element of the space 
$W^{1,2,2}(0,T;H^{s-1/2}(\partial \Omega),H^{s-5/2}(\partial \Omega))$, 
with $s\geq 2$. In this case, we can obtain an extension in the space 
$W^{1,2,2}(0,T;H^s(\Omega),H^{s-2}(\Omega))$ (cf. Theorem 3.2 of 
\cite{fursikov1}) and, if we want to ensure that $W^{1,2,2}(0,T;H^s(\Omega),
H^{s-2}(\Omega))$ $\subset W$, we can take, for instance, $s\geq 3$.
\end{rem}

\subsection{Part 2: Galerkin approximation}

In this part we will construct a sequence of approximations that will converge to a solution of problem 
(\ref{eq:Asystem12}). So, let $\{\omega_n\}_{n \in \mathbb{N}} \subset \widetilde{X}$ be 
a dense subset of independent vectors of $\widetilde{X}$, {\color{red} such 
that ${\omega_n}_{|_{\Gamma_R}} \in L^{\infty}(\Gamma_R)$, $\forall n \in \mathbb{N}$}, 
which we can assume orthonormal in $L^2(\Omega)$. {\color{red} We also assume} that 
the projection
\begin{equation}
P_N(v)=\sum_{k=1}^N \left(\int_{\Omega} \omega_k v \, d \mathbf{x} \right) 
\omega_k
\end{equation}
is selfadjoint and 
$\|{P_N}_{|_Z}\|_{\mathcal{L}(Z,Z)} \leq 1$, $\forall N \in \mathbb{N}$, 
where $Z$ is a Banach space, as given in Lemma \ref{lemma1n}. 
Then, for $N \in \mathbb{N}$, we denote by:
\begin{equation}
\xi_N=\sum_{n=1}^N \xi_n^N(t) \omega_n,
\end{equation}
where the coefficients $\xi_n^N(t)$, $n=1,\ldots,N$, are such that $\xi_N$ is the solution 
of the following differential equation:
\begin{equation} \label{eq:Asystem3n}
\hspace{-0.3cm}\begin{array}{r}
\displaystyle \int_{\Omega} \frac{\partial \xi_N}{\partial t} \omega_k \, d \mathbf{x} +
\int_{\Omega} \mathbf{v} \cdot \nabla \xi_N \omega_k \, d \mathbf{x} + 
K \int_{\Omega} \nabla \xi_N \cdot \nabla \omega_k \, d \mathbf{x} +  
b_1^A \int_{\Gamma_A} \xi_N \omega_k \, d \gamma \\  
\displaystyle 
+b_1^R \int_{\Gamma_R} \xi_N \omega_k \, d \gamma
+b_2^R \int_{\Gamma_R} |\xi_N+\zeta_{D}|^3(\xi_N+\zeta_{D}) 
\omega_k \, d \gamma = 
\int_{\Omega} H_{D}\,  \omega_k \, d \mathbf{x}
\\  
\displaystyle
+b_1^A \int_{\Gamma_A}   g_{D}^{A}\, \omega_k \, d \gamma
+b_1^R \int_{\Gamma_R}   g_{D}^{R}\, \omega_k \, d \gamma
+b_2^R \int_{\Gamma_R} \varphi^4 \omega_k \, d \gamma, \;
\forall k=1,\ldots,N,
\end{array}
\end{equation}
which can be rewritten in the following standard formulation:
\begin{equation} \label{eq:Asystem4}
\left\{\begin{array}{l}
\displaystyle \frac{d\mathbf{y}}{dt}=\mathbf{F}(\mathbf{y}(t),t), \; a.e.\, t \in ]0,T[, \\  
\displaystyle \mathbf{y}(0)=\mathbf{y}_0,
\end{array}\right.
\end{equation}
where:
\begin{equation}
\mathbf{y}(t)=\left(\xi_1^N(t),\, \xi_2^N(t),\ldots,\xi_N^N(t) \right)^T,
\end{equation}
\begin{equation}
\mathbf{y}_0=\left((\xi_0^N,\omega_1),(\xi_0^N,\omega_2),\ldots, (\xi_0^N,\omega_N) \right)^T,
\end{equation}
\begin{equation}
\mathbf{F}(\mathbf{y},t)=\left(\begin{array}{c}
\langle f(t),\omega_1 \rangle -a(t;\mathbf{y} \cdot \boldsymbol{\omega},\omega_1) \\ 
\langle f(t),\omega_2 \rangle -a(t;\mathbf{y} \cdot \boldsymbol{\omega},\omega_2) \\ 
\vdots \\ 
\langle f(t),\omega_N \rangle -a(t;\mathbf{y} \cdot \boldsymbol{\omega},\omega_N)
\end{array}\right),
\end{equation}
\begin{equation}
\boldsymbol{\omega}=\left( 
\omega_1,\omega_2,\ldots,\omega_N
\right)^T,
\end{equation}
\begin{equation}
\begin{array}{rcl}
\displaystyle
\langle f(t),\omega_k \rangle &&=\displaystyle
\int_{\Omega} H_{D}(t)\,  \omega_k \, d \mathbf{x}
+b_1^A \int_{\Gamma_A}   g_{D}^{A}(t)\, \omega_k \, d \gamma
\\  
+ && \displaystyle
b_1^R \int_{\Gamma_R}   g_{D}^{R}(t)\, \omega_k \, d \gamma
+b_2^R \int_{\Gamma_R} \varphi^4(t) \omega_k \, d \gamma, 
\end{array}
\end{equation}
\begin{equation} 
\begin{array}{rcl}
\displaystyle
a(t;\mathbf{y} \cdot \boldsymbol{\omega},\omega_k)&&=
\displaystyle
\int_{\Omega} \mathbf{v}(t) \cdot \nabla (\mathbf{y} \cdot \boldsymbol{\omega}) \omega_k \, d \mathbf{x} \\ 
\displaystyle &&+ 
K \int_{\Omega} \nabla (\mathbf{y} \cdot \boldsymbol{\omega} ) \cdot \nabla \omega_k \, d \mathbf{x} 
\\  
+ && \displaystyle
b_1^A \int_{\Gamma_A} (\mathbf{y} \cdot \boldsymbol{\omega}) \omega_k \, d \gamma 
+b_1^R \int_{\Gamma_R} (\mathbf{y} \cdot \boldsymbol{\omega}) \omega_k \, d \gamma
\\  
+ && \displaystyle
b_2^R \int_{\Gamma_R} |(\mathbf{y} \cdot \boldsymbol{\omega})+
\zeta_{D}(t)|^3((\mathbf{y} \cdot \boldsymbol{\omega})+\zeta_{D}(t)) 
\omega_k \, d \gamma,
\end{array}
\end{equation}
where $k=1,\ldots,N$, and $\xi_0^N=P_N(\theta^0-\zeta_{D}(0))$ is the projection 
of $\theta^0-\zeta_{D}(0)$ onto $\widetilde{X}_N= \langle \{\omega_1,\ldots,\omega_N\} \rangle $. 

Thus, we say than an element $\xi_N \in W^{1,5/4}(0,T;\widetilde{X}_N)$ is a solution 
of system (\ref{eq:Asystem3n}) if it satisfies the ordinary differential equation problem (\ref{eq:Asystem4}). 
We must recall here that $\|\xi_0^N\|_{L^2(\Omega)} \leq \|\theta^0-\zeta_{D}(0)\|_{L^2(\Omega)}$, 
$\forall N \in \mathbb{N}$, and that, if $\xi_N$ is a solution 
of problem (\ref{eq:Asystem3n}), then $\mathbf{y} \in \mathcal{C}([0,T];\mathbb{R}^N)$.

\begin{lem} Within the framework established in Hypothesis \ref{hypo1n}, 
there exits a constant $C>0$ independent of $N$ such that:
\begin{equation} \label{eq:bound1}
\begin{array}{r}
\displaystyle 
\| \xi_N \|^2_{L^{\infty}(0,T;L^2(\Omega))}+
\|\xi_N\|^2_{L^2(0,T;H^1(\Omega))} +
\|\xi_N\|^5_{L^5(0,T;L^5(\Gamma_S))}  \\  
\displaystyle \leq
C\Big[ \|\theta^0-\zeta_{D}(0)\|_{L^2(\Omega)}^2+
\|H_D\|^2_{L^2(0,T;L^2(\Omega))}
+ \|g_{D}^{A}\|_{L^2(0,T;L^2(\Gamma_A))}^2
\\  
\displaystyle
+
\|g_{D}^{R}\|_{L^2(0,T;L^2(\Gamma_R)}^2+
\|\zeta_{D}\|_{L^5(0,T;L^5(\Gamma_R))}^5
+\|\varphi\|_{L^5(0,T;L^5(\Gamma_R))}^5 \Big],
\end{array}
\end{equation}
\begin{equation} \label{eq:bound2} \hspace{-.7cm}
\begin{array}{r}
\displaystyle 
\left\|
\frac{d \xi_N}{dt}
\right\|_{L^{5/4}(0,T;Z')} \leq C \Big[
\|H_D\|_{L^2(0,T;L^2(\Omega))}
+\|g_{D}^{A}  \|_{L^2(0,T;L^2(\Gamma_A))}
\\  
\displaystyle
+\| g_{D}^{R}  \|_{L^2(0,T;L^2(\Gamma_R))}
+\|\varphi\|^4_{L^5(0,T;L^5(\Gamma_R))}
+\|\mathbf{v}\|_{L^{10/3}(0,T;[L^3(\Omega)]^3)}
\|\xi_N\|_{L^2(0,T;\widetilde{X})}
\\  
\displaystyle +
\|\xi_N\|_{L^2(0,T;\widetilde{X})}
+\|\xi_N\|^4_{L^5(0,T;L^5(\Gamma_R))}
+\|\zeta_{D}\|^4_{L^5(0,T;L^5(\Gamma_R))}\Big].
\end{array}
\end{equation}
\end{lem}

\begin{pf}  Multiplying
(\ref{eq:Asystem3n}) by $\xi_k^N(t)$, summing in $k$ and adding 
to both sides the term
\begin{equation}
b_2^R \int_{\Gamma_R} 
|\xi_N(s)+\zeta_{D}(s)|^3(\xi_N(s)+\zeta_{D}(s)) \zeta_{D}(s) \, d \mathbf{x},
\end{equation}
we have:
\begin{equation}
\begin{array}{r}
\displaystyle 
 \frac{1}{2}\frac{d}{dt} \|\xi_N(s)\|^2_{L^2(\Omega)}
+K \| \nabla \xi_N(s)\|^2_{[L^2(\Omega)]^3} 
+b_2^R \| \xi_N(s)+\zeta_{D}(s)\|^5_{L^5(\Gamma_R)}
\\ 
\displaystyle
\leq 
\|H_{D}(s) \|_{L^2(\Omega)}\|\xi_N(s)\|_{L^2(\Omega)}
+b_1^A \|g_{D}^{A}(s) \|_{L^2(\Gamma_A)} \| \xi_N(s) \|_{L^2(\Gamma_A)}
\\ 
\displaystyle
+b_1^R \|g_{D}^{R}(s) \|_{L^2(\Gamma_R)} \| \xi_N(s) \|_{L^2(\Gamma_R)}
+b_2^R \|\varphi(s)\|^4_{L^5(\Gamma_R)} \|\xi_N(s)\|_{L^5(\Gamma_R)} 
\\ 
\displaystyle
+b_2^R \|\xi_N(s)+\zeta_{D}(s)\|^4_{L^5(\Gamma_R)}
\|\zeta_{D}(s)\|_{L^5(\Gamma_R)}.
\end{array}
\end{equation}
As a consequence of the continuity of trace operator and of the inequalities of Young and Poincare, we obtain:
\begin{equation} \nonumber
\begin{array}{rcl}
\displaystyle 
\|H_{D}(s) \|_{L^2(\Omega)}\|\xi_N(s)\|_{L^2(\Omega)} 
&\leq & 
\displaystyle
\frac{C_1}{\epsilon_1} \|H_{D}(s) \|_{L^2(\Omega)}^2+ \epsilon_1 
\| \nabla \xi_N(s)\|^2_{[L^2(\Omega)]^3} ,
 \vspace{0.1cm}\\ 
\displaystyle 
\| g_{D}^{A}(s)  \|_{L^2(\Gamma_A)} \| \xi_N(s)\|_{L^2(\Gamma_A)}
&\leq &
\displaystyle
\frac{C_2}{\epsilon_2} \| g_{D}^{A}(s)  \|_{L^2(\Gamma_A)}^2 +
\epsilon_2 \| \nabla \xi_N(s)\|^2_{[L^2(\Omega)]^3} ,
 \vspace{0.1cm}\\ 
\displaystyle 
\| g_{D}^{R}(s)  \|_{L^2(\Gamma_R)} \| \xi_N(s)\|_{L^2(\Gamma_R)}
&\leq & 
\displaystyle
\frac{C_3}{\epsilon_3}\| g_{D}^{R}(s)  \|_{L^2(\Gamma_R)}^2+
\epsilon_3 \| \nabla \xi_N(s)\|^2_{[L^2(\Omega)]^3} ,
 \vspace{0.1cm}\\ 
\displaystyle 
\|\varphi(s)\|^4_{L^5(\Gamma_R)} \|\xi_N(s)\|_{L^5(\Gamma_R)}
&\leq & 
\displaystyle
\frac{C_4}{\sqrt[4]{\epsilon_4}} 
\|\varphi(s)\|^5_{L^5(\Gamma_R)}  +
\epsilon_4 \|\xi_N(s)\|_{L^5(\Gamma_R)}^5,
 \vspace{0.1cm}\\ 
\displaystyle 
\|\xi_N(s)+\zeta_{D}(s)\|^4_{L^5(\Gamma_R)}
\|\zeta_{D}(s)\|_{L^5(\Gamma_R)}
&\leq & 
\displaystyle
\frac{C_5}{\sqrt[4]{\epsilon_5}} \|\zeta_{D}(s)\|_{L^5(\Gamma_R)}^5
\\ && \displaystyle +\epsilon_5 
 \|\xi_N(s)+\zeta_{D}(s)\|^5_{L^5(\Gamma_R)},
\end{array}
\end{equation}
where $\epsilon_k$, $k=1,\ldots,5$, are arbitrary strictly positive numbers, and 
$C_k$, $k=1,\ldots,5$, are constants that may depend on the trace operator and 
Young and Poincare's inequalities. If we take $\epsilon_k$, $k=1,2,3$, 
such that $K-\sum_{k=1}^3 \epsilon_k=K/2$ and $\epsilon_5=b_2^R/2$, then (renaming the constants if necessary): 

\begin{equation}
\begin{array}{r}
\displaystyle 
\frac{1}{2}\frac{d}{dt} \|\xi_N(s)\|^2_{L^2(\Omega)}
+\frac{K}{2} \| \nabla \xi_N(s)\|^2_{[L^2(\Omega)]^3} 
+\frac{b_2^R}{2} \| \xi_N(s)+\zeta_{D}(s)\|^5_{L^5(\Gamma_R)}
\\ 
\displaystyle
 \leq 
C_1 \|H_{D}(s) \|_{L^2(\Omega)}^2
+C_2 \| g_{D}^{A}(s) \|_{L^2(\Gamma_A)}^2
+C_3 \|g_{D}^{R}(s)  \|_{L^2(\Gamma_R)}^2
\\ 
\displaystyle
+\frac{C_4}{\sqrt[4]{\epsilon_4}} \|\varphi(s)\|^5_{L^5(\Gamma_R)} 
+\epsilon_4 \|\xi_N(s)\|_{L^5(\Gamma_R)}^5
+C_5 \|\zeta_{D}(s)\|_{L^5(\Gamma_R)}^5.
\end{array}
\end{equation}
Adding to both sides $\frac{b^S_2}{2}\|\zeta_{D}(s)\|_{L^5(\Gamma_R)}^5$, 
using the inequality $(a+b)^5 \leq 16 a^5+ 16 b^5$ (for $a,b\geq 0$), and taking 
$\epsilon_5=b^R_2/64$: 
\begin{equation}
\begin{array}{r}
\displaystyle 
\frac{1}{2}\frac{d}{dt} \|\xi_N(s)\|^2_{L^2(\Omega)}
+\frac{K}{2} \| \nabla \xi_N(s)\|^2_{[L^2(\Omega)]^3} 
+\frac{b_2^R}{64} \| \xi_N(s)\|^5_{L^5(\Gamma_R)}
\\ 
 \displaystyle
 \leq 
 C_1 \|H_{D}(s) \|_{L^2(\Omega)}^2 + 
 C_2 \| g_{D}^{A}(s)  \|_{L^2(\Gamma_R)}^2+
C_3 \|g_{D}^{R}(s)  \|_{L^2(\Gamma_R)}^2
\\ 
\displaystyle+
C_4 \|\varphi(s)\|^5_{L^5(\Gamma_R)} + 
C_5 \|\zeta_{D}(s)\|_{L^5(\Gamma_R)}^5.
\end{array}
\end{equation}
Finally, integrating over the time interval $[0,t]$ (renaming again the constants):
\begin{equation} \label{eq:bound3}
\begin{array}{r}
\displaystyle 
\|\xi_N(t)\|^2_{L^2(\Omega)} + 
\int_0^t \|\nabla \xi_N(s)\|^2_{[L^2(\Omega)]^3} \, ds + 
\int_0^t \| \xi_N(s)\|^5_{L^5(\Gamma_R)}\, ds
\\ 
 \displaystyle
 \leq 
C \Big[ 
\|\theta^0-\zeta_{D}(0)\|^2_{L^2(\Omega)}
+\int_0^t \|H_{D}(s) \|_{L^2(\Omega)}^2 \, ds
+\int_0^t \| g_{D}^{A}(s)  \|_{L^2(\Gamma_R)}^2 \, ds
\\ 
 \displaystyle
+\int_0^t \|\varphi(s)\|^5_{L^5(\Gamma_R)} \, ds 
+\int_0^t \|\zeta_{D}(s)\|_{L^5(\Gamma_R)}^5\, ds
\Big],
\end{array}
\end{equation}
where we have used that $\|\xi_0^N\|_{L^2(\Omega)} \leq \|\theta^0-\zeta_{D}(0)\|_{L^2(\Omega)}$, 
$\forall N \in \mathbb{N}$. Finally, applying Gronwall's Lemma, we obtain that there 
exists a positive constant 
$C$ independent of $N$ such that (\ref{eq:bound1}) is satisfied. 

For obtaining (\ref{eq:bound2}) it is sufficient to apply Holder inequality and bear in mind
the fact that the projection operator $P_N$ onto $\widetilde{X}_N$ 
is bounded independently of $N$:
\begin{equation}  \nonumber
\begin{array}{r}
\displaystyle 
\left<\frac{d \xi_N}{dt}, v \right>=
\left<P_N\left(\frac{d \xi_N}{dt}\right), v \right>=
\left<\frac{d \xi_N}{dt}, P_N (v) \right>  \\
\displaystyle
\leq 
C \Big[
\|H_D\|_{L^2(0,T;L^2(\Omega))}
+\|g_{D}^{A}  \|_{L^2(0,T;L^2(\Gamma_A))}
\\  
\displaystyle
+\| g_{D}^{R}  \|_{L^2(0,T;L^2(\Gamma_R))}
+\|\varphi\|^4_{L^5(0,T;L^5(\Gamma_R))}
\\  
\displaystyle
+\|\mathbf{v}\|_{L^{10/3}(0,T;[L^3(\Omega)]^3)}
\|\xi_N\|_{L^2(0,T;\widetilde{X})} 
\\  
\displaystyle
+ \|\xi_N\|_{L^2(0,T;\widetilde{X})}+\|\xi_N\|^4_{L^5(0,T;L^5(\Gamma_R))}
\\  
\displaystyle
+\|\zeta_{D}\|^4_{L^5(0,T;L^5(\Gamma_R))}\Big] 
\|P_N\|_{\mathcal{L}(Z,Z)} \|v\|_{L^5(0,T;Z)}.
\end{array}
\end{equation}
\hfill $\blacksquare$
\end{pf}

\subsection{Part 3: Existence of solution for the Galerkin approximation}

Now we will demonstrate that there exists, for each $N \in \mathbb{N}$, a unique absolutely continuous 
solution $\xi_N$ of equation (\ref{eq:Asystem3n}). 

\begin{lem}\label{lemma3n} Within the framework established in Hypothesis 
\ref{hypo1n}, there exists a unique absolutely continuous solution defined on the 
whole time interval $[0,T]$ of Cauchy problem (\ref{eq:Asystem4}).
\end{lem}

\begin{pf} To state this lemma we can apply the Caratheodory theorem for ordinary differential
equations (cf., for example, Theorem 5.2 of \cite{hale1}). Indeed, $\mathbf{F}(\cdot,t)$ is continuous for any $t \in [0,T]$, and 
$\mathbf{F}(\mathbf{y}, \cdot ) \in L^{5/4}(0,T)$ for any $\mathbf{y} \in \mathbb{R}^N$. Then, 
given an open ball $B$ in $\mathbb{R}^N$, if we prove that there exist two functions 
$m_B,\, l_B \in L^1(0,T)$ such that:
\begin{equation} \label{eq:caratheocond} 
\begin{array}{rcl}
\displaystyle 
\|\mathbf{F}(\mathbf{y},t)\| &\leq & \displaystyle 
m_B(t),\; a.e. \; t \in ]0,T[, \; \forall\, \mathbf{y} \in B, \\ 
\displaystyle
\|\mathbf{F}(\mathbf{y}_1,t)-
\mathbf{F}(\mathbf{y}_2,t)\| & \leq & 
\displaystyle l_B(t) \|\mathbf{y}_1-\mathbf{y}_2\|,\; a.e. \; t \in ]0,T[, \; 
\forall\, \mathbf{y}_1,\, \mathbf{y}_2 \in B,
\end{array}
\end{equation}
we can conclude that problem (\ref{eq:Asystem4}) has a unique absolutely continuous 
solution, which can be extended to the boundary of $]0,T[ \times B$. It is worthwhile mentioning here that, if 
\begin{equation} 
\begin{array}{r}
diam(B)^2> C\Big[ \|\theta^0-\zeta_{D}(0)\|^2_{L^2(\Omega)}+
\|H_D\|^2_{L^2(0,T;L^2(\Omega))}
 \\  
\displaystyle
+ \|g_{D}^{A}\|^2_{L^2(0,T;L^2(\Gamma_A))}
+ \|g_{D}^{R}\|^2_{L^2(0,T;L^2(\Gamma_R)}
 \\  
\displaystyle
+ \|\zeta_{D}\|_{L^5(0,T;L^5(\Gamma_R))}^5
+\|\varphi\|_{L^5(0,T;L^5(\Gamma_R))}^5 \Big],
\end{array}
\end{equation}
then the solution $\mathbf{y}$ cannot reach the boundary of $B$ because of the 
\textit{a priori} estimate (\ref{eq:bound1}) and the fact that $\|\xi_N(t)\|^2_{L^2(\Omega)}=
\|\mathbf{y}(t)\|^2$. The first of the Caratheodory conditions (\ref{eq:caratheocond}) can be obtained by 
applying Holder inequality (in fact, we obtain that $m_B \in L^{5/4}(0,T)$). 
In the other hand, for the second Caratheodory condition, we can use the following 
inequality (straightforward consequence of the mean value Theorem):
\begin{equation} \label{eq:trick1}
b|b|^3-a|a|^3=4|c|^3(b-a),
\end{equation}
for $c=\lambda b + (1-\lambda) a$ and $\lambda \in (0,1)$. So, we 
obtain the following inequality:

\begin{equation} \nonumber 
\begin{array}{r}
\displaystyle
|F_k(\mathbf{y}_1,t)-
F_k(\mathbf{y}_2,t)| = \left|
a(t;\mathbf{y}_1 \cdot  \boldsymbol{\omega},\omega_k)-
a(t;\mathbf{y}_2 \cdot  \boldsymbol{\omega},\omega_k)
\right|
\\  
 \displaystyle
\leq
\displaystyle C_1\Big[
\|\mathbf{v}(t)\|_{[L^3(\Omega)]^3} \|\nabla \boldsymbol{\omega}\|_{[L^2(\Omega)]^{N \times 3}} 
\|\omega_k\|_{L^6(\Omega)}
\\  
\displaystyle
+ \|\nabla \boldsymbol{\omega}\|_{[L^2(\Omega)]^{N \times 3}} 
\|\nabla \omega_k\|_{L^2(\Omega)}
+\|\boldsymbol{\omega}\|_{[L^2(\Gamma_A)]^N} \|\omega_k\|_{L^2(\Gamma_A)}
\\  
\displaystyle
+\|\boldsymbol{\omega}\|_{[L^2(\Gamma_R)]^N} \|\omega_k\|_{L^2(\Gamma_R)} \Big]
\|\mathbf{y}_1-\mathbf{y}_2\|
\\  
 \displaystyle
+C_2 \Big[ 
\Big(\max\{\|\mathbf{y}_1\|,\|\mathbf{y}_2\|\}^3
\|\boldsymbol{\omega}\|_{[L^5(\Gamma_R)]^N}^3
\\  
\displaystyle
+\|\xi_D(t)\|_{L^5(\Gamma_R)}^3
\Big) 
\|\boldsymbol{\omega}\|_{[L^5(\Gamma_R)]^N}
\|\omega_k\|_{L^5(\Gamma_R)}
\Big]\|\mathbf{y}_1-\mathbf{y}_2\|.
\end{array}
\end{equation}
Therefore, we can conclude the existence of a 
function $l_B \in L^1(0,T)$ (in fact, $l_B \in L^{5/3}(0,T))$ such that the second 
of the Caratheodory conditions (\ref{eq:caratheocond}) is achieved. \hfill $\blacksquare$
\end{pf}

\subsection{Part 4: Convergence of the Galerkin approximation}

In previous subsections we have seen that there exists a bounded 
sequence $\{\xi_N\}_{N \in \mathbb{N}} \subset \widetilde{W}$ of 
solutions of problem (\ref{eq:Asystem4}). In this subsection we will 
pass to the limit and obtain a solution of equation (\ref{eq:Asystem2}). 

\begin{lem} \label{lemma4} There exists a subsequence of 
$\{\xi_N\}_{N \in \mathbb{N}}$, still denoted in the same way, such that:
\begin{enumerate}
\item $\displaystyle \xi_N \rightharpoonup \xi$ in $L^2(0,T;\widetilde{X})$, 
\item $\displaystyle \xi_N \rightharpoonup^* \xi$ in $L^{\infty}(0,T;L^2(\Omega))$, 
\item $\displaystyle \xi_N \rightarrow \xi$ in $L^{10/3-\epsilon}(0,T;L^{10/3-\epsilon}(\Omega))$,
\item $\displaystyle \xi_N \rightarrow \xi$ in $L^2(0,T;L^2(\Gamma_A))$,
\item {\color{red} $\displaystyle \xi_N \rightarrow \xi$ in $L^4(0,T;L^4(\Gamma_R))$}.
\end{enumerate}
\end{lem}

\begin{pf} Thanks to the boundedness in $\widetilde{W}$ of the sequence 
$\{\xi_N\}_{N \in \mathbb{N}}$, we obtain the first two convergences. 
The third and fourth limits are a direct consequence of Aubin, and 
Lions Lemma (cf. Lemma 7.7 of \cite{Roubicek1}) and the 
compactness of $H^1(\Omega)$, respectively, in 
$L^{6-\epsilon}(\Omega)$ and $L^{4-\epsilon}(\partial \Omega)$, $\forall \epsilon>0$. 
Finally, the fifth convergence is a consequence of Lemma \ref{lemma1n}.  \hfill $\blacksquare$
\end{pf}

\begin{lem} \label{lemma5} If $\{\xi_N\}_{N \in \mathbb{N}}$ is a bounded sequence in 
$\widetilde{W}$, then there exists a subsequence of 
$\{\xi_N\}_{N \in \mathbb{N}}$, still denoted in the same way, such that, for all 
{\color{red} $ \phi \in L^{\infty}(0,T;L^{\infty}(\Gamma_R))$}:
\begin{equation} 
\begin{array}{c}
\displaystyle
\lim_{N \to \infty} \int_0^T \int_{\Gamma_R} |\xi_N+\zeta_{D}|^3(\xi_N+\zeta_{D}) 
\phi\, d \gamma\, dt \\  
\displaystyle
= \int_0^T \int_{\Gamma_R} |\xi+\zeta_{D}|^3(\xi+\zeta_{D}) 
\phi\, d \gamma \, dt.
\end{array}
\end{equation}
\end{lem}

\begin{pf} Using the same technique that we have employed in the proof of 
Lemma \ref{lemma3n}, from the {\color{red} strong convergence of $\{\xi_N\}_{N \in \mathbb{N}}$ 
to $\xi$ in $L^4(0,T;L^4(\Gamma_R))$ we have}:
{\color{red}\begin{equation}
\begin{array}{r}
\displaystyle 
\lim_{N \to \infty} 
\int_0^T \int_{\Gamma_R}\big[ |\xi_N+\zeta_{D}|^3(\xi_N+\zeta_{D}) -
 |\xi+\zeta_{D}|^3(\xi+\zeta_{D})  \big] \phi\, d \gamma \, dt 
\\  \displaystyle
\leq \lim_{N \to \infty} C \int_0^T \int_{\Gamma_R} 
\big[|\xi_N-\xi|^4+|\xi_D-\xi|^3|\xi_N-\xi| \big] \phi\, d \gamma \, dt \\
\displaystyle 
\leq C \lim_{n \to \infty} \big[\|\xi_N-\xi\|^4_{L^4(0,T;L^4(\Gamma_R))}\\ \displaystyle 
+
\|\xi_D-\xi\|^3_{L^4(0,T;L^4(\Gamma_R))} \|\xi_N-\xi\|_{L^4(0,T;L^4(\Gamma_R))}
 \big] \|\phi\|_{L^{\infty}(0,T;L^{\infty}(\Gamma_R))}=0.
\end{array}
\end{equation}}
\hfill $\blacksquare$
\end{pf}

\subsection{Part 5: Proof of the main result of the Appendix}

Now, we can demonstrate the Theorem \ref{Atheo1}: 

\begin{pf} We will divide the proof into three 
parts, in the first part we will pass to the limit in the Galerkin 
approximation in order to obtain a solution 
for the system (\ref{eq:Asystem2}), in the second part, we will derive 
the estimates (\ref{eq:boundcont1}) and (\ref{eq:boundcont2}) and, 
finally, in the third part we will prove the uniqueness of solution.

First, for a fixed index $k \in \mathbb{N}$, if we multiply (\ref{eq:Asystem3n}) by a scalar function $\psi$ 
continuously differentiable on $[0,T]$, such that $\psi(T)=0$, integrate with respect to $t$, and 
integrate by parts, we have, $\forall N \geq k$:
\begin{equation} \nonumber 
\begin{array}{r}
\displaystyle -\int_0^T \int_{\Omega} 
\xi_N(t) \frac{d \psi}{dt}(t) \omega_k \, d \mathbf{x} \, dt
+\int_0^T \int_{\Omega} \mathbf{v}(t) \cdot \nabla \xi_N(t) \psi(t) \omega_k \, d \mathbf{x}\, dt 
\vspace{0.1cm}\\  
\displaystyle
+K \int_0^T\int_{\Omega} \nabla \xi_N(t) \cdot \nabla \omega_k \psi(t) \, d \mathbf{x}\, dt 
+b_1^A \int_0^T \int_{\Gamma_A} \xi_N(t) \psi(t) \omega_k \, d \gamma \, dt
\vspace{0.1cm}\\  
\displaystyle 
+b_1^R \int_0^T \int_{\Gamma_R} \xi_N(t) \psi(t) \omega_k \, d \gamma \, dt
\\ \displaystyle
+b_2^R \int_0^T \int_{\Gamma_R} |\xi_N(t)+\zeta_{D}(t)|^3(\xi_N(t)+\zeta_{D}(t)) 
\psi(t)\omega_k \, d \gamma\, dt 
\vspace{0.1cm}\\  
\displaystyle
= \int_{\Omega} \xi_0^N \omega_k \psi(0)\, d \mathbf{x}
+\int_0^T \int_{\Omega} H_{D}(t)\,  \psi(t) \omega_k \, d \mathbf{x}\, dt
+b_1^A \int_0^T \int_{\Gamma_A}   g_{D}^{A}(t)\, \psi(t) \omega_k \, d \gamma\, dt
\vspace{0.1cm}\\  
\displaystyle
+b_1^R\int_0^T \int_{\Gamma_R}   g_{D}^{R}(t)\,\psi(t) \omega_k \, d \gamma\, dt
+b_2^R \int_0^T \int_{\Gamma_R} \varphi^4(t) \psi(t)\omega_k \, d \gamma\, dt.
\end{array}
\end{equation}
The passage to the limit for $N \to \infty$ in the integrals of the left-hand side is due to the Lemmas 
\ref{lemma4} and \ref{lemma5}. We observe also that $\xi_0^N \to \theta^0-\zeta_D(0)$ in $L^2(\Omega)$. 
Hence, we find in the limit:

\begin{equation} \label{eq:Asystem5} 
\begin{array}{r}
\displaystyle -\int_0^T \int_{\Omega} 
\xi(t) \frac{d \psi}{dt}(t) \eta \, d \mathbf{x} \, dt
+\int_0^T \int_{\Omega} \mathbf{v}(t) \cdot \nabla \xi(t) \psi(t) \eta \, d \mathbf{x}\, dt 
\vspace{0.1cm}\\  
\displaystyle
+K \int_0^T\int_{\Omega} \nabla \xi(t) \cdot \nabla \omega_k \psi(t) \, d \mathbf{x}\, dt 
+b_1^A \int_0^T \int_{\Gamma_A} \xi(t) \psi(t) \eta \, d \gamma \, dt
\vspace{0.1cm}\\  
\displaystyle 
+b_1^R \int_0^T \int_{\Gamma_R} \xi(t) \psi(t) \eta \, d \gamma \, dt \\
\displaystyle
+b_2^R \int_0^T \int_{\Gamma_R} |\xi(t)+\zeta_{D}(t)|^3(\xi(t)+\zeta_{D}(t)) 
\psi(t)\eta \, d \gamma\, dt 
\vspace{0.1cm}\\  
\displaystyle
= \int_{\Omega} (\theta^0-\zeta_D(0)) \eta \psi(0)\, d \mathbf{x}
+\int_0^T \int_{\Omega} H_{D}(t)\,  \psi(t) \eta \, d \mathbf{x}\, dt
\vspace{0.1cm}\\  
\displaystyle
+b_1^A \int_0^T \int_{\Gamma_A}   g_{D}^{A}(t)\, \psi(t) \eta \, d \gamma\, dt
+b_1^R\int_0^T \int_{\Gamma_R}   g_{D}^{R}(t)\,\psi(t) \eta \, d \gamma\, dt
\vspace{0.1cm}\\  
\displaystyle
+b_2^R \int_0^T \int_{\Gamma_R} \varphi^4(t) \psi(t)\eta \, d \gamma\, dt,
\end{array}
\end{equation}
for each $ \eta \in \widetilde{X}$ which is a finite lineal combination of elements $\omega_k$. Since each term 
of above expression depends linearly and continuously on $\eta$, for the norm of $\widetilde{X}$, 
previous equality remains still valid, by continuity, for each $\eta \in \widetilde{X}$. Now, writing in particular 
(\ref{eq:Asystem5}) for $\psi=\phi \in \mathcal{D}(0,T)$, we obtain the variational formulation (\ref{eq:Asystem2}). 
Finally, we can prove that $\xi(0)=\theta^0-\zeta_D(0)$ multiplying (\ref{eq:Asystem2}) by the same $\psi$ as before, 
integrating by parts with respect to $t$, and comparing with (\ref{eq:Asystem5}).

Then, multiplying inequality (\ref{eq:bound3}) by $\phi$, with $\phi \in \mathcal{D}(0,T)$, $\phi\geq 0$, 
and integrating in $[0,T]$ we have:
\begin{equation} \label{eq:bound4}  
\begin{array}{r}
\displaystyle 
\int_0^T \Big\{
\|\xi_N(t)\|^2_{L^2(\Omega)} + 
\int_0^t \|\nabla \xi_N(s)\|^2_{[L^2(\Omega)]^3} \, ds + 
 \vspace{0.1cm}\\ 
 \displaystyle
\int_0^t \| \xi_N(s)\|^5_{L^5(\Gamma_R)}\, ds \Big\} \phi(t) \, dt
 \leq 
C  \int_0^T\Big\{
\|\theta^0-\zeta_{D}(0)\|^2_{L^2(\Omega)}
 \vspace{0.1cm}\\ 
 \displaystyle
+\int_0^t \|H_{D}(s) \|_{L^2(\Omega)}^2 \, ds
+\int_0^t \| g_{D}^{A}(s)  \|_{L^2(\Gamma_R)}^2 \, ds
     \vspace{0.1cm}\\ 
 \displaystyle
+\int_0^t \|\varphi(s)\|^5_{L^5(\Gamma_R)} \, ds 
+\int_0^t \|\zeta_{D}(s)\|_{L^5(\Gamma_R)}^5\, ds\Big\}
\phi(t) \,dt,
\end{array}
\end{equation}
from which, taking into account that the norm of a reflexive Banach 
space is weakly lower semicontinuous, we can pass to the inferior limit thanks to convergences of 
Lemma \ref{lemma4}:
\begin{equation} \label{eq:boundcont3} 
\begin{array}{r}
\displaystyle 
\int_0^T \Big\{
\|\xi(t)\|^2_{L^2(\Omega)} + 
\int_0^t \|\nabla \xi(s)\|^2_{[L^2(\Omega)]^3} \, ds + 
 \vspace{0.1cm}\\ 
 \displaystyle
\int_0^t \| \xi(s)\|^5_{L^5(\Gamma_R)}\, ds \Big\} \phi(t) \, dt

 \leq 
C  \int_0^T\Big\{
\|\theta^0-\zeta_{D}(0)\|^2_{L^2(\Omega)}
 \vspace{0.1cm}\\ 
 \displaystyle
+\int_0^t \|H_{D}(s) \|_{L^2(\Omega)}^2 \, ds
+\int_0^t \| g_{D}^{A}(s)  \|_{L^2(\Gamma_R)}^2 \, ds
+\int_0^t \|\varphi(s)\|^5_{L^5(\Gamma_R)} \, ds 
 \vspace{0.1cm}\\ 
 \displaystyle
+\int_0^t \|\zeta_{D}(s)\|_{L^5(\Gamma_R)}^5\, ds\Big\}
\phi(t) \,dt, \quad \forall \phi \in \mathcal{D}(0,T),\ \phi\geq 0.
\end{array}
\end{equation}
Thus, we obtain the following energy inequality for a.e. $t \in (0,T)$:
\begin{equation} \label{eq:boundcont4}  
\begin{array}{r}
\displaystyle 
\|\xi(t)\|^2_{L^2(\Omega)} + 
\int_0^t \|\nabla \xi(s)\|^2_{[L^2(\Omega)]^3} \, ds + 
\int_0^t \| \xi(s)\|^5_{L^5(\Gamma_R)}\, ds 
    \vspace{0.1cm}\\ 
 \displaystyle
 \leq 
C  \Big[\int_0^T
\|\theta^0-\zeta_{D}(0)\|^2_{L^2(\Omega)}
+\int_0^t \|H_{D}(s) \|_{L^2(\Omega)}^2 \, ds
     \vspace{0.1cm}\\ 
 \displaystyle
+\int_0^t \| g_{D}^{A}(s)  \|_{L^2(\Gamma_R)}^2 \, ds
+\int_0^t \|\varphi(s)\|^5_{L^5(\Gamma_R)} \, ds 
+\int_0^t \|\zeta_{D}(s)\|_{L^5(\Gamma_R)}^5\, ds\Big].
\end{array}
\end{equation}
Finally, (\ref{eq:boundcont1}) can be derived from above expression thanks to the Gronwall's Lemma, and 
estimate (\ref{eq:boundcont2}) is a direct consequence of Holder inequality.  

Now, we will prove the uniqueness of solution. Let us assume the existence 
of two solutions $\xi_1$ and $\xi_2$ for problem 
(\ref{eq:Asystem2}), and define $\xi_{12}=\xi_1-\xi_2$. We have that $\xi_{12}\in \widetilde{W}$, 
$\xi_{12}(0)=0$, a.e. $\mathbf{x} \in \Omega$, and that satisfies the following variational formulation:
\begin{equation}
\begin{array}{r}
\displaystyle \int_{\Omega} \frac{\partial \xi_{12}}{\partial t} \eta \, d \mathbf{x} +
\int_{\Omega} \mathbf{v} \cdot \nabla \xi_{12} \eta \, d \mathbf{x} 
\\ \displaystyle
+ 
K \int_{\Omega} \nabla \xi_{12} \cdot \nabla \eta \, d \mathbf{x} +  
b_1^A \int_{\Gamma_A} \xi_{12} \eta \, d \gamma  
+b_1^R \int_{\Gamma_R} \xi_{12} \eta \, d \gamma
\\ \displaystyle 
+b_2^R \int_{\Gamma_R} \big[|\xi_1+\zeta_{D}|^3(\xi_1+\zeta_{D}) 
-|\xi_2+\zeta_{D}|^3(\xi_2+\zeta_{D}) \big] \eta \, d \gamma 
\\  
\displaystyle 
= 0,  \quad a.e.\; t \in ]0,T[, \quad \forall \eta \in \widetilde{X}.
\end{array}
\end{equation}

From a direct computation, we observe that $\frac{d \xi_{12}}{dt}=\alpha+\beta+\gamma$, where $\alpha \in 
L^{2}(0,T;\widetilde{X}')$, $\beta \in L^{5/4}(0,T;L^{5/4}(\Gamma))$ 
and $\gamma \in L^{4/3}(0,T;L^{5/6}(\Omega))$, defined by:
\begin{equation} \nonumber
\langle \alpha,z \rangle =
-K \int_0^T \int_{\Omega}  \nabla \xi_{12} \cdot \nabla z \, d \mathbf{x} \, dt
-b_1^A \int_0^T \int_{\Gamma_A} \xi_{12} z \, d \gamma\, dt
-b_1^R \int_0^T \int_{\Gamma_R} \xi_{12} z \, d \gamma\, dt,
\end{equation} 
$\forall z \in \widetilde{X}$, $\ \ \beta=
-b_2^R \big[|{\xi_1}_{|_{\Gamma_R}}
+{\zeta_{D}}_{|_{\Gamma_R}}|^3({\xi_1}_{|_{\Gamma_R}}+
{\zeta_{D}}_{|_{\Gamma_R}}) 
-|{\xi_2}_{|_{\Gamma_R}}+{\zeta_{D}}_{|_{\Gamma_R}}|^3(
{\xi_2}_{|_{\Gamma_R}}+{\zeta_{D}}_{|_{\Gamma_R}}) \big] $, and 
$\gamma=-\mathbf{v} \cdot \nabla \xi_{12}$. 

Now, for each $\delta>0$, we define the following function:
\begin{equation}
\kappa_{\delta}(r)=\left\{ \begin{array}{ll}
1 & \mbox{ if } r\geq \delta, \\  
\frac{r}{\delta} & \mbox{ if } -\delta\leq r \leq \delta,\\ 
-1 & \mbox{ if }  r \leq -\delta,
\end{array}\right.
\end{equation}
and its primitive,
\begin{equation} \label{Kd}
K_{\delta}(r)=\int_0^r \kappa_{\delta}(s) \, ds = 
\left\{\begin{array}{ll}
r-\frac{\delta}{2} & \mbox{ if } r \geq \delta, \\ 
\frac{r^2}{2 \delta} & \mbox{ if } -\delta \leq r \leq \delta, \\ 
-r-\frac{\delta}{2} & \mbox{ if } r \leq -\delta.
\end{array} \right.
\end{equation}
We can extend the results proved in Corollary 9 of \cite{fran2} to our case 
(we just need to add the term in $L^{5/4}(0,T;L^{5/4}(\Gamma_R))$ in 
the proof of Lemma 6 of \cite{fran2} and then we can use this result to prove 
Corollary 9) and then, we have:
\begin{equation} \hspace{-.9cm}
\begin{array}{r}
\displaystyle 
\int_{\Omega} K_{\delta}(\xi_{12}(t)) \, dx -\int_{\Omega} 
K_{\delta}(\xi_{12}(0)) \, d \mathbf{x} =  
\int_0^t \langle \alpha(s),\kappa_{\delta}(\xi_{12}(s)) \rangle_{\widetilde{X}',
\widetilde{X}}\, ds
\\ 
\displaystyle 
+\int_0^t \int_{\Gamma_R} \beta(s) \kappa_{\delta}(\xi_{12}(s)) \, 
d \gamma \, ds
+
\int_0^t \int_{\Omega} \gamma(s) \kappa_{\delta}(\xi_{12}(s)) \, d \mathbf{x} \, ds
, \quad a.e.\, t \in ]0,T[.
\end{array}
\end{equation}
By a direct evaluation of previous expression, we have:
\begin{equation}  \nonumber
\begin{array}{r}
\displaystyle 
\int_{\Omega} K_{\delta}(\xi_{12}(t))\, d \mathbf{x}
+K \int_0^t \int_{\Omega} \nabla \xi_{12}(s) \cdot \nabla \kappa_{\delta}(\xi_{12}(s)) \, d \mathbf{x} \, ds 
\\ \displaystyle 
+\int_0^t  \int_{\Omega} \mathbf{v}(s) \cdot \nabla \xi_{12}(s) \kappa_{\delta}(\xi_{12}(s)) \, d \mathbf{x}\, dt
+ b_1^A \int_0^t \int_{\Gamma_A} \xi_{12}(s) \kappa_{\delta}(\xi_{12}(s))\, d \gamma \, ds\\ \displaystyle
+b_1^R \int_{0}^t \int_{\Gamma_R} \xi_{12}(s) \kappa_{\delta}(\xi_{12}(s))\, d \gamma \, ds
+b_2^R \int_0^t \int_{\Gamma_R} \big[|\xi_1(s)+\zeta_{D}(s)|^3(\xi_1(s)+\zeta_{D}(s))\\ 
\displaystyle  
-|\xi_2(s)+\zeta_{D}(s)|^3(\xi_2(s)+\zeta_{D}(s)) \big] \kappa_{\delta}(\xi_{12}(s)) \, d \gamma\, ds =0,
\end{array}
\end{equation}
a.e. $t \in ]0,T[$. Moreover, it is obvious that: 
\begin{equation}  \hspace{-.5cm} \nonumber
\begin{array}{r}
\displaystyle
\int_{\Omega} \mathbf{v}(s) \cdot \nabla \xi_{12}(s) \kappa_{\delta}(\xi_{12}(s)) \, d \mathbf{x}=0, \; a.e.\; s \in ]0,T[,
 \vspace{0.1cm}\\  
\displaystyle
\int_{\Omega} \nabla \xi_{12}(s) \cdot \nabla \kappa_{\delta}(\xi_{12}(s)) \, d \mathbf{x}=
\int_{\Omega} \kappa_{\delta}'(\xi_{12}(s)) \|\nabla \xi_{12}(s)\|^2 \, d \mathbf{x} \geq 0, \; a.e.\; s \in ]0,T[,
 \vspace{0.1cm}\\ 
\displaystyle 
\int_{\Gamma_A} \xi_{12}(s) \kappa_{\delta}(\xi_{12}(s))\, d \gamma \geq 0,\;\; a.e.\; s \in ]0,T[,
 \vspace{0.1cm}\\ 
\displaystyle 
\int_{\Gamma_R} \xi_{12}(s) \kappa_{\delta}(\xi_{12}(s))\, d \gamma \geq 0,\;\; a.e.\; s \in ]0,T[,
 \vspace{0.1cm}\\ 
\displaystyle
\int_{\Gamma_R} \big[|\xi_1(s)+\zeta_{D}(s)|^3(\xi_1(s)+\zeta_{D}(s)) 
-|\xi_2(s)+\zeta_{D}(s)|^3(\xi_2(s)+\zeta_{D}(s)) \big] \kappa_{\delta}(\xi_{12}(s)) \, d \gamma
\\ 
\displaystyle
=4\int_{\Omega}|c(s)|^3 \xi_{12}(s) \kappa_{\delta}(\xi_{12}(s)) \, d \gamma \geq 0,\;\; a.e.\; s \in ]0,T[,
\end{array}
\end{equation}
for $c(s)=\lambda(s)\xi_1(s)+(1-\lambda(s))\xi_2(s)$ and $\lambda(s) \in (0,1) $, a.e. $s\in ]0,T[$. Thus, we can deduce 
that:
\begin{equation}
\int_{\Omega} K_{\delta}(\xi_{12}(t))\, d \mathbf{x} \leq 0,\; a.e.\; t \in ]0,T[.
\end{equation}
Finally, from the definition (\ref{Kd}) of $K_{\delta}$ it is straightforward that $0\leq |r|-K_{\delta}(r)\leq \frac{\delta}{2}$, 
$\forall r \in \mathbb{R}$, and then:
\begin{equation}
\|\xi_{12}(t)\|_{L^1(\Omega)} \leq \frac{\delta}{2} \, | \Omega |
+\int_{\Omega} K_{\delta}(\xi_{12}(t))\, d \mathbf{x}, \quad \forall \delta > 0,
\end{equation}
with $ | \Omega |$ denoting the volume of $ \Omega $, 
which implies that $\xi_{12}(\mathbf{x},t)=0$, a.e. $(\mathbf{x},t) \in \Omega \times ]0,T[$,
and, consequently, $\xi_{1}=\xi_{2}$.  \hfill $\blacksquare$
\end{pf}

%%-----------------------------
%%      your bibliography
%%-----------------------------

\end{document}